\documentstyle[amscd,amssymb,verbatim,epsf]{amsart}

\begin{document}
\theoremstyle{plain}
\newtheorem{Thm}{Theorem}
\newtheorem{Cor}{Corollary}
\newtheorem{Con}{Conjecture}
\newtheorem{Main}{Main Theorem}
\newtheorem{Lem}{Lemma}
\newtheorem{Prop}{Proposition}

\theoremstyle{definition}
\newtheorem{Def}{Definition}
\newtheorem{Note}{Note}

\theoremstyle{remark}
\newtheorem{notation}{Notation}
\renewcommand{\thenotation}{}

\errorcontextlines=0
\numberwithin{equation}{section}
\renewcommand{\rm}{\normalshape}%

\title[On the space of oriented geodesics of hyperbolic 3-space]%
   {On the space of oriented geodesics of hyperbolic 3-space}
\author{Nikos Georgiou}
\address{Nikos Georgiou\\
          Department of Computing and Mathematics\\
          Institute of Technology, Tralee\\
          Clash\
          Tralee\\
          Co. Kerry\\
          Ireland.}
\email{nikos.georgiou@@research.ittralee.ie}
\author{Brendan Guilfoyle}
\address{Brendan Guilfoyle\\
          Department of Computing and Mathematics\\
          Institute of Technology, Tralee\\
          Clash\\
          Tralee\\
          Co. Kerry\\
          Ireland.}
\email{brendan.guilfoyle@@ittralee.ie}

\keywords{Kaehler structure, hyperbolic 3-space, isometry group}
\subjclass{Primary: 51M09; Secondary: 51M30}
\date{9th February, 2007}

\begin{abstract}
We construct a K\"ahler structure (${\mathbb{J}},\Omega,{\mathbb{G}}$) on the space ${\mathbb{L}}({\mathbb{H}}^3)$ of oriented geodesics of hyperbolic 3-space ${\mathbb{H}}^3$ and investigate its properties. We prove that (${\mathbb{L}}({\mathbb{H}}^3),{\mathbb{J}})$ is biholomorphic to ${\mathbb{P}}^1\times{\mathbb{P}}^1-\overline{\Delta}$, where $\overline{\Delta}$ is the reflected diagonal, and that the K\"ahler metric ${\mathbb{G}}$ is of neutral signature, conformally flat and scalar flat. We establish that the identity component of the isometry group of the metric ${\mathbb{G}}$ on ${\mathbb{L}}({\mathbb{H}}^3)$ is isomorphic to the identity component of the hyperbolic isometry group. Finally, we show that the geodesics of ${\mathbb{G}}$ correspond to ruled minimal surfaces in ${\mathbb{H}}^3$, which are totally geodesic iff the geodesics are null.
\end{abstract}

\maketitle

The space ${\mathbb{L}}({\mathbb{M}}^3)$ of oriented geodesics on a 3-manifold ${\mathbb{M}}^3$ of constant curvature is a 
4-dimensional manifold which carries a natural complex structure ${\mathbb{J}}$. In the case where ${\mathbb{M}}^3$ is 
Euclidean 3-space ${\mathbb{E}}^3$, this complex structure can be traced back to Weierstrass \cite{weier} and 
Whittaker \cite{whitt}, with its modern re-emergence occurring in Hitchen's study of monopoles on ${\mathbb{E}}^3$ 
\cite{hitch}.

More recently, this structure has been supplemented by a compatible symplectic structure, so that 
${\mathbb{L}}({\mathbb{M}}^3)$ inherits a natural K\"ahler structure. This has been investigated when 
${\mathbb{M}}^3={\mathbb{E}}^3$ and ${\mathbb{M}}^3={\mathbb{E}}^3_1$ \cite{gak4} \cite{gak5} \cite{gak6} and the
purpose of this paper is to study the hyperbolic 3-space case ${\mathbb{M}}^3={\mathbb{H}}^3$.

From a topological point of view ${\mathbb{L}}({\mathbb{M}}^3)$ is homeomorphic to S$^2\times$S$^2-\Delta$, where
$\Delta$ is the diagonal. However, from holomorphic point of view we show that:

\vspace{0.1in}

\noindent{\bf Theorem 1}:

{\it
The complex surface (${\mathbb{L}}({\mathbb{H}}^3),{\mathbb{J}})$ is biholomorphic to 
${\mathbb{P}}^1\times{\mathbb{P}}^1-\overline{\Delta}$, where $\overline{\Delta}$ is the reflected diagonal 
(see Definition \ref{d:refdiag}).
}

\vspace{0.1in}

The ${\mathbb{P}}^1$ 
in the Theorem refers to the boundary of the Poincar\'e ball model of ${\mathbb{H}}^3$, considered as the 
past and future infinities of the oriented geodesics, from which ${\mathbb{L}}({\mathbb{H}}^3)$
inherits its complex structure.

We then turn to the K\"ahler metric ${\mathbb{G}}$ and prove:

\vspace{0.1in}

\noindent{\bf Theorem 2}:

{\it
The K\"ahler metric ${\mathbb{G}}$ is of neutral signature, conformally flat and scalar flat.
}

\vspace{0.1in}

We also show that, despite the ($++--$) signature, this metric on ${\mathbb{L}}({\mathbb{H}}^3)$ faithfully reflects the 
hyperbolic metric $g$ on ${\mathbb{H}}^3$, in the following sense:

\newpage

\noindent{\bf Theorem 3}:

{\it
The identity component of the isometry group of the metric ${\mathbb{G}}$ on ${\mathbb{L}}({\mathbb{H}}^3)$ 
is isomorphic to the identity component of the hyperbolic isometry group.  
}

\vspace{0.1in}

A curve in ${\mathbb{L}}({\mathbb{H}}^3)$ is a 1-parameter family of oriented geodesics in ${\mathbb{H}}^3$: a ruled 
surface. Our final result characterises the ruled surfaces that arise as geodesics in ${\mathbb{L}}({\mathbb{H}}^3)$:

\vspace{0.1in}
\noindent{\bf Theorem 4}:

{\it
The geodesics of the K\"ahler metric ${\mathbb{G}}$ are generated by the 1-parameter subgroups of the isometry group
of  ${\mathbb{G}}$. 

A ruled surface generated by a geodesic of ${\mathbb{G}}$ is a minimal surface in ${\mathbb{H}}^3$, and the
geodesic is null iff the ruled surface is totally geodesic.
}

\vspace{0.1in}

In the next section we describe the space of oriented geodesics of hyperbolic 3-space from a topological and a
differentiable point of view, using the ball and upper half-space models of ${\mathbb{H}}^3$.
In section 2 we define and investigate the K\"ahler structure on ${\mathbb{L}}({\mathbb{H}}^3)$ and prove Theorems 1 
and 2. The proof of Theorem 3, which is contained in section 3, consists of a number of steps, formulated as propositions. 
We first find the Killing vectors of ${\mathbb{G}}$. We then compute the action induced on ${\mathbb{L}}({\mathbb{H}}^3)$
by isometries of ${\mathbb{H}}^3$, and prove that the infinitesimal generators of this action coincide precisely 
with the Killing vectors of ${\mathbb{G}}$. Finally, we study the geodesics of the neutral K\"ahler metric 
and the ruled surfaces they generate in ${\mathbb{H}}^3$ in section 4. 

The uniqueness of this K\"ahler structure has recently been established by Salvai \cite{salvai2} ({\it cf}. \cite{salvai1} 
for the ${\mathbb{E}}^n$ case). Indeed, a number of our results overlap with those of Salvai, who utilises 
techniques of Lie groups to obtain his results. Our approach is particularly geared to the study of surfaces in 
${\mathbb{L}}({\mathbb{H}}^3)$ and, given recent interest in flat and CMC surfaces in ${\mathbb{H}}^3$ \cite{kok} 
\cite{roit} \cite{small}, we hope that this line of inquiry will prove fruitful. 

\vspace{0.2in}

\section{The space of oriented geodesics of ${\mathbb{H}}^3$}

\begin{Def}
Let ${\mathbb{L}}({\mathbb{H}}^3)$ be the space of oriented geodesics in ${\mathbb{H}}^3$.
\end{Def}

The topology of ${\mathbb{L}}({\mathbb{H}}^3)$ is most easily seen using the Poincar\'e ball model. This has underlying
space
\[
 B^3=\{(y^{1},y^{2},y^{3})\in {\mathbb{R}}^3\;|\; (y^{1})^2+(y^{2})^2+(y^{3})^2< 1\},
\]
for standard coordinates $(y^{1},y^{2},y^{3})$ on ${\mathbb{R}}^3$, with hyperbolic metric
\[
ds^2=\frac{4[(dy^1)^{2}+(dy^2)^{2}+(dy^3)^{2}]}{[1-(y^{1})^2-(y^{2})^2-(y^{3})^2]^2}.
\]

\begin{Def}\label{d:refdiag}
Let $\iota:$S$^2\rightarrow$S$^2$ be the antipodal map and define the {\it reflected diagonal} by
\[
\overline{\Delta}=\{(\mu_1,\mu_2)\in S^2\times S^2\;|\;\mu_1=\iota(\mu_2)\;\}.
\]
\end{Def}

We can now identify the space of oriented geodesics of ${\mathbb{H}}^3$:

\vspace{0.2in}

\begin{Prop}
The space ${\mathbb{L}}({\mathbb{H}}^3)$ of oriented geodesics on hyperbolic 3-space is homeomorphic to S$^2\times$S$^2-\overline{\Delta}$.
\end{Prop}
\begin{pf}
Consider the unit ball model of ${\mathbb{H}}^3$. In this model, the geodesics are either diameters, or circles which are asymptotically 
orthogonal to the 
boundary 2-sphere. An oriented geodesic can thus be uniquely identified by its beginning and end point on the boundary. Moreover, any 
ordered pair of points on the boundary 2-sphere define a unique oriented geodesic, as long as the points are distinct. Thus, the space of 
oriented geodesics is homeomorphic to S$^2\times$S$^2-\{{\mbox{diag}}\}$. 

In fact, for geometric
reasons which will become clear below, we will identify an oriented geodesic by the direction of its tangent vector at past and future infinity - see Figure 1. 
Since these directions are inward and outward pointing (respectively), we see that the oriented geodesics
can also be identified with S$^2\times$S$^2$ minus antipodal directions, as claimed.
\end{pf}

\vspace{0.1in}
\setlength{\epsfxsize}{4.5in}
\begin{center}
   \mbox{\epsfbox{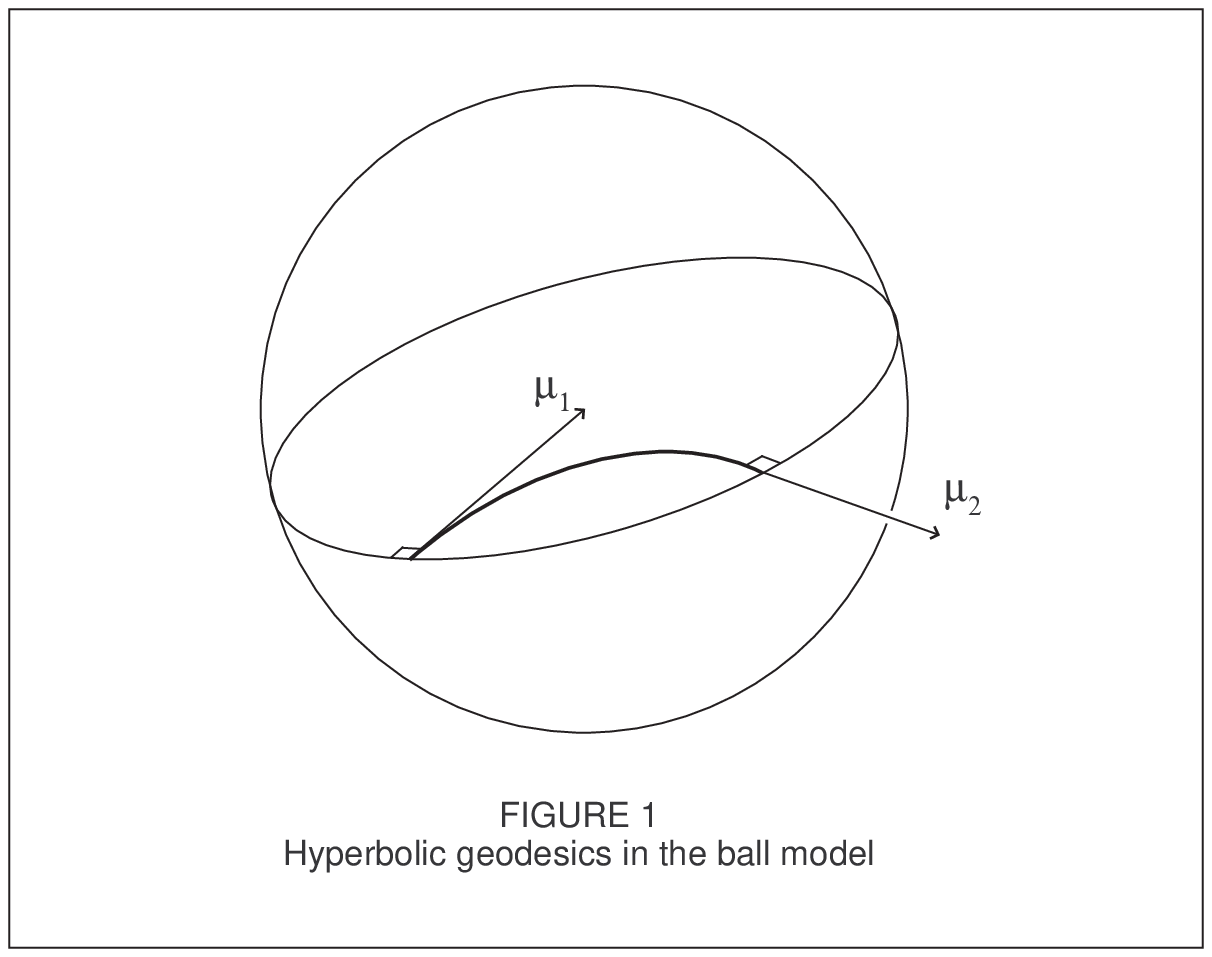}}
\end{center}
\vspace{0.1in}

For computational purposes we use the upper half-space model of ${\mathbb{H}}^3$.  Thus, the underlying space is 
\[
{\mathbb{R}}^3_{+}=\{\;(x^{0},x^{1},x^{2})\in {\mathbb{R}}^3 \;|\; x^{0}>0\;\},
\]
for standard coordinates $(x^{0},x^{1},x^{2})$ on ${\mathbb{R}}^3$. In these coordinates the hyperbolic metric $g$ 
has expression:
\[
d\tilde{s}^2=\frac{(dx^{0})^2+(dx^{1})^2+(dx^{2})^2}{(x^{0})^2}.
\]

This is related to the ball model by the mapping $P\colon {\mathbb{R}}^3_{+}\rightarrow  B^3\colon (x^{0},x^{1},x^{2}) \mapsto (y^{1},y^{2},y^{3})$
defined by
\[
y^{1}=\frac{2x^1}{(x^0+1)^2+(x^1)^2+(x^2)^2},\quad 
y^{2}=\frac{2x^2}{(x^0+1)^2+(x^1)^2+(x^2)^2},
\]
\begin{equation}\label{e:ballcoords}
y^{3}=\frac{(x^0)^2+(x^1)^2+(x^2)^2-1}{(x^0+1)^2+(x^1)^2+(x^2)^2}.
\end{equation}
The map $P$ is a diffeomorphism. In fact, it is an isometry: $P^{*}(ds^2)=d\tilde{s}^2$ \cite{Rat}.

We now describe the geodesics of ${\mathbb{H}}^3$ in this model.

\vspace{0.2in}

\begin{Prop}\label{p:geo1}
The geodesics of ${\mathbb{H}}^3$ that are not parallel to the $x^{0}-$axis are
\[
x^{0}=\frac{1}{\sqrt{c_{1}^{2}+c_{2}^{2}}\cosh(r+r_0)}, \qquad  x^{1}=\frac{c_{1}}{c_{1}^{2}+c_{2}^{2}}\tanh (r+r_0)+c_{3},
\]
\[
x^{2}=\frac{c_{2}}{c_{1}^{2}+c_{2}^{2}}\tanh(r+r_0)+c_{4},
\]
where $r$ is the arc length of the geodesic, $c_{1}, c_{2}, c_{3},c_{4},r_0\in{\mathbb{R}}$ and $c_{1},c_{2}$ are not both $0$. Geodesics parallel to the $x^{0}-$axis
are given by
\[
x^{0}(r)=e^{r+r_{0}} \qquad x^{1}(r)=c_{3} \qquad x^{2}(r)=c_{4}, 
\]
for $c_{3}, c_{4}, r_0\in{\mathbb{R}}$.
\end{Prop}
\begin{pf}
Let $(x^{0},x^{1},x^{2})$ be local coordinates on the space ${\mathbb{H}}^{3}$ with metric $g$ as defined above.
The only non-vanishing christoffel symbols of the metric are
\[
\Gamma_{00}^{0}=\Gamma_{01}^{1}=\Gamma_{02}^{2}=-\frac{1}{x^{0}},\qquad  \Gamma_{11}^{0}=\Gamma_{22}^{0}=\frac{1}{x^{0}}.
\]
The geodesic equations (using the summation convention here and throughout)
\[
\frac{d^2 x^{k}}{dr^2}+\Gamma_{ij}^{k}\frac{dx^{i}}{dr}\frac{dx^{j}}{dr}=0,\quad k=0,1,2 ,
\]
then turn out to be:
\[
\frac{d^2 x^{0}}{dr^2}+\frac{1}{x^{0}}\left[\left(\frac{dx^{1}}{dr}\right)^2+\left(\frac{dx^{2}}{dr}\right)^2-\left(\frac{dx^{0}}{dr}\right)^2\right]=0,
\]
\[ 
\frac{d^2 x^{1}}{dr^2}-\frac{2}{x^{0}}\frac{dx^{1}}{dr}\frac{dx^{0}}{dr}=0, \qquad \frac{d^2 x^{2}}{dr^2}-\frac{2}{x^{0}}\frac{dx^{2}}{dr}\frac{dx^{0}}{dr}=0 ,
\]
where $r$ is an affine parameter along the geodesic.

These can be integrated to yield the first integrals:
\[
I_{1}=\frac{(\dot{x}^{0})^2+(\dot{x}^{1})^2+(\dot{x}^{2})^2}{(x^{0})^2},\qquad I_{2}=\frac{2\dot{x}^{2}}{(x^{0})^2},\qquad  I_{3}=\frac{2\dot{x}^{1}}{(x^{0})^2} ,
\]
where the dot denotes differentiation with respect to $r$. Thus, $I_{1},I_{2},I_{3}$ are constant along any
geodesic.

By parameterising the geodesic by arc-length we can set $I_{1}=1$. Let $I_{2}=2c_{1}$ and $I_{3}=2c_{2}$, so that 
\begin{equation}\label{e:rel1}
(\dot{x}^{0})^2+(\dot{x}^{1})^2+(\dot{x}^{2})^2=(x^{0})^2,\qquad \dot{x}^{1}=c_{1}(x^{0})^2,\qquad \dot{x}^{2}=c_{2}(x^{0})^2.
\end{equation}
Combining these equations, we then get 
\[
c_{1}^2 (x^{0})^4+c_{2}^2 (x^{0})^4+(\dot{x}^{0})^2=(x^{0})^2,
\]
or, rearranging:
\[
(\dot{x}^{0})^2=(x^{0})^2-(c_{1}^2+c_{2}^2) (x^{0})^4.
\]
Integrating, we have that:
\begin{equation}\label{e:int1}
\int{\frac{dx^{0}}{x^{0}\sqrt{1-K(x^{0})^2}}}=r+r_{0} \qquad K=c_{1}^2+c_{2}^2.
\end{equation}
For $c_{1}\neq 0$ or $c_{2}\neq 0$ we find that
\[
\int{\frac{dx^{0}}{x^{0}\sqrt{1-K(x^{0})^2}}}=\frac{1}{2}\log\left|\frac{1-\sqrt{1-K(x^{0})^2}}{1+\sqrt{1-K(x^{0})^2}}\right|=\frac{1}{2}\log\left(\frac{1-\sqrt{1-K(x^{0})^2}}{1+\sqrt{1-K(x^{0})^2}}\right).
\]
So we have
\[
1-\sqrt{1-K(x^{0})^2}=e^{2(r+r_{0})}(1+\sqrt{1-K(x^{0})^2}),
\]
and hence
\[
x^{0}=\frac{1}{\sqrt{c_{1}^2+c_{2}^2}\cosh(r+r_{0})},
\]
as claimed.

Now, from the second of equation (\ref{e:rel1})
\[
\dot{x}^{1}=c_{1}(x^{0})^2,
\]
we obtain that
\[
\dot{x}^{1}=\frac{c_{1}}{(c_{1}^2+c_{2}^2)\cosh^{2}(r+r_{0})}.
\]
Similarly we obtain
\[
\dot{x}^{2}=\frac{c_{2}}{(c_{1}^2+c_{2}^2)\cosh^{2}(r+r_{0})},
\]
and integrating we finally get
\[
x^{1}(r)=\frac{c_{1}}{c_{1}^2+c_{2}^2}\tanh(r+r_{0})+c_{3},\qquad x^{2}(r)=\frac{c_{2}}{c_{1}^2+c_{2}^2}\tanh(r+r_{0})+c_{4},
\]
as claimed.

The case $c_{1}=c_{2}=0$ follows easily by integration of (\ref{e:int1}).
\end{pf}

\vspace{0.2in}

We see that the geodesics in ${\mathbb{H}}^3$, where $c_{1},c_{2}$ are not both 0, are semi-circles in 3-space  
${\mathbb{R}}^3$ with centre $(0,c_{3},c_{4})$ and radius $(c_{1}^2+c_{2}^2)^{\scriptstyle{-\frac{1}{2}}}$. If we 
let 
\[
\xi=c_{1}+ic_{2}, \qquad \eta=c_{3}+ic_{4},
\]
then the geodesics can be labelled as shown in Figure 2.

\vspace{0.1in}
\setlength{\epsfxsize}{4.5in}
\begin{center}
   \mbox{\epsfbox{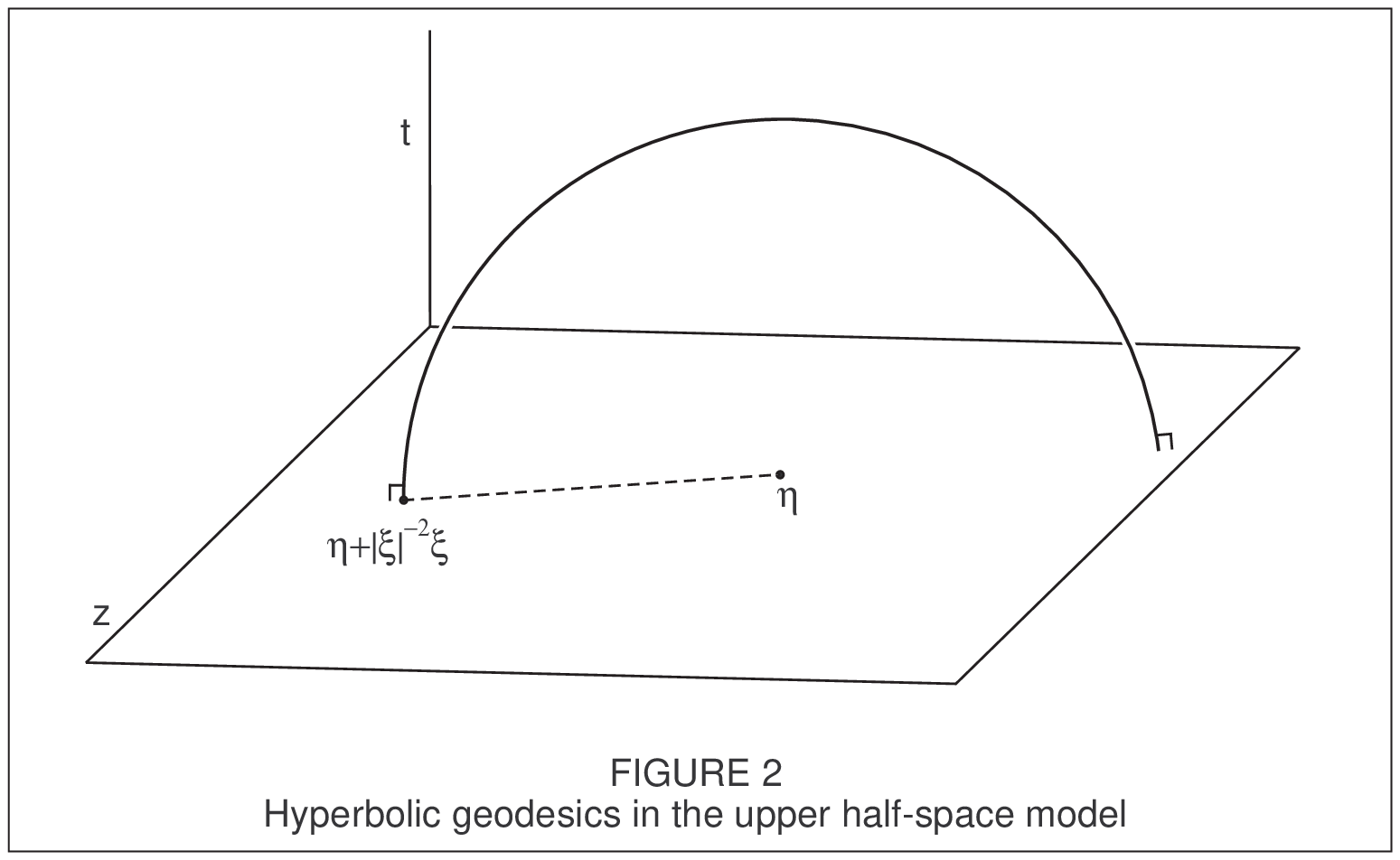}}
\end{center}
\vspace{0.1in}

In fact, initially we will define a K\"ahler structure only on an open subset of ${\mathbb{L}}({\mathbb{H}}^3)$. We 
define this subset as follows.

\begin{Def}
Let ${\mathbb{U}}\subset{\mathbb{L}}({\mathbb{H}}^3)$ be the set of oriented geodesics in the upper-half space model that
are not parallel to the $x^0-$axis. 

\end{Def}

Introducing complex coordinates $t=x^{0}, z=x^{1}+ix^{2}$, so that
\[
\frac{\partial}{\partial z}=\frac{1}{2}\left(\frac{\partial}{\partial x^{1}}-i\frac{\partial}{\partial x^{2}}\right), \quad \frac{\partial}{\partial t}=\frac{\partial}{\partial x^{0}},
\]
the metric tensor $g$ becomes
\[
ds^2=\frac{dzd\bar{z}+dt^2}{t^2}.
\]

From Proposition \ref{p:geo1} the points of ${\mathbb{U}}$ are the oriented geodesics in ${\mathbb{H}}^3$ that
are given in terms of the constants $\xi$ and $\eta$ by
\begin{equation}\label{e:phi}
t=\frac{1}{|\xi|\cosh (r+r_0)} \qquad z=\eta + \frac{\tanh (r+r_0)}{\bar{\xi}},
\end{equation}
for $\xi\in{\mathbb{C}}-\{0\},\eta\in{\mathbb{C}}$ and $r_0\in{\mathbb{R}}$. 

Thus, ($\xi,\eta$) are local coordinates on ${\mathbb{U}}\subset{\mathbb{L}}({\mathbb{H}}^3)$ and by 
setting $r_0=0$ we fix the parameterisation on these geodesics. In the next section we define a K\"ahler 
structure on ${\mathbb{U}}$ and extend it to all of ${\mathbb{L}}({\mathbb{H}}^3)$. First, we must
explicitly identify the tangent space T${\mathbb{U}}$ with the orthogonal Jacobi fields along the associated geodesics in 
${\mathbb{H}}^3$ \cite{kan}.

\begin{Def}
Let $\Phi :{\mathbb{U}}\times {\mathbb{R}}\rightarrow {\mathbb{H}}^3$ be the map $\Phi(\xi,\eta,r)=\left(t(\xi,\eta,r),z(\xi,\eta,r)\right)$ given by equations (\ref{e:phi}) with $r_0=0$.
\end{Def}

For later use we note that:

\begin{Prop}\label{p:dphi}
The derivative $D\Phi:T_{(\xi,\eta,r)}({\mathbb{U}}\times {\mathbb{R}})\rightarrow T_{\Phi(\xi,\eta,r)}({\mathbb{H}}^3)$ is given by
\begin{align}
D\Phi_{(\xi,\eta,r)}\left(\frac{\partial}{\partial \xi}\right)&=-\frac{\tanh r}{\xi^2}\frac{\partial}{\partial \bar{z}}-\frac{1}{2\xi|\xi|\cosh r}\frac{\partial}{\partial t},\nonumber \\
D\Phi_{(\xi,\eta,r)}\left(\frac{\partial}{\partial \eta}\right)&=\frac{\partial}{\partial z}, \nonumber \\
D\Phi_{(\xi,\eta,r)}\left(\frac{\partial}{\partial r}\right)&=\frac{1}{\bar{\xi}\cosh^2 r}\frac{\partial}{\partial z}+\frac{1}{\xi\cosh^2 r}\frac{\partial}{\partial \bar{z}}-\frac{\sinh r}{|\xi|\cosh^2 r}\frac{\partial}{\partial t}. \nonumber 
\end{align}
\end{Prop}
\begin{pf}
This is found by differentiation of $\Phi$.
\end{pf}

\begin{Def}
A {\it null frame} at a point $p$ in ${\mathbb{H}}^3$ is a trio of vectors
$e_{(0)}$, $e_{(+)}$, $e_{(-)}\in{\mathbb{C}}\otimes \mbox{T}_p{\mathbb{H}}^3$ such that:
\[
e_{(0)}=\overline{e_{(0)}}, \qquad\qquad e_{(+)}=\overline{e_{(-)}},
\]
\[
g(e_{(0)},e_{(0)})=g(e_{(+)},e_{(-)})=1, \qquad\qquad g(e_{(0)},e_{(+)})=0,
\] 
where the hyperbolic inner product $g$ is extended bilinearly
over ${\mathbb{C}}$.

Given an oriented geodesic in ${\mathbb{H}}^3$, an {\it adapted} null frame is a null frame along the geodesic such 
that $e_{(0)}$ is the tangent to the geodesic.
\end{Def}

\begin{Prop}\label{p:adapt}
An adapted null frame to the oriented geodesic $(\xi,\eta)\in{\mathbb{U}}$ is given by
\[
e_{(0)}=\frac{1}{\bar{\xi}\cosh^2 r}\frac{\partial}{\partial z}+\frac{1}{\xi\cosh^2 r}\frac{\partial}{\partial \bar{z}}-\frac{\sinh r}{|\xi|\cosh^2 r}\frac{\partial}{\partial t} ,
\]
\[
e_{(+)}=\frac{1}{\sqrt{2}\cosh^2 r}\left[-\frac{e^{-r}}{\bar{\xi}}\frac{\partial}{\partial z}+\frac{e^{r}}{\xi}\frac{\partial}{\partial \bar{z}}\right]+\frac{1}{\sqrt{2}|\xi|\cosh^2 r}\frac{\partial}{\partial t} .
\]
\end{Prop}
\begin{pf}
We note that $e_{(0)}=D\Phi_{(\xi,\eta,r)}\left(\frac{\partial}{\partial r}\right)$ and a straightforward computation gives
\[
g(e_{(\alpha)},e_{(\beta)})=\begin{bmatrix} 1 & 0 & 0 \\ 0 & 0 & 1 \\ 0 & 1 & 0 \end{bmatrix} \qquad \alpha,\beta=0,+,-.
\]
Thus we have an adapted null frame as claimed.
\end{pf}
\noindent{\bf Note}: An orthonormal frame along the geodesic $(\xi,\eta)$ is given by
\begin{align}
e_{(0)}&=-\frac{1}{\sqrt{c_{1}^{2}+c_{2}^{2}}}\frac{\sinh r}{\cosh^2 r}\frac{\partial}{\partial x^{0}}+\frac{c_{1}}{c_{1}^{2}+c_{2}^{2}}\frac{1}{\cosh^2 r}\frac{\partial}{\partial x^{1}}+\frac{c_{2}}{c_{1}^{2}+c_{2}^{2}}\frac{1}{\cosh^2 r}\frac{\partial}{\partial x^{2}},\nonumber \\
e_{(1)}&=\frac{1}{\sqrt{c_{1}^{2}+c_{2}^{2}}}\frac{1}{\cosh^2 r}\frac{\partial}{\partial x^{0}}+\frac{c_{1}}{c_{1}^{2}+c_{2}^{2}}\frac{\sinh r}{\cosh^2 r}\frac{\partial}{\partial x^{1}}+\frac{c_{2}}{c_{1}^{2}+c_{2}^{2}}\frac{\sinh r}{\cosh^2 r}\frac{\partial}{\partial x^{2}}, \nonumber \\
e_{(2)}&=-\frac{c_{2}}{c_{1}^{2}+c_{2}^{2}}\frac{1}{\cosh r}\frac{\partial}{\partial x^{1}}+\frac{c_{1}}{c_{1}^{2}+c_{2}^{2}}\frac{1}{\cosh r}\frac{\partial}{\partial x^{2}}. \nonumber
\end{align}
In particular,
\[
e_{(+)}=\frac{1}{\sqrt{2}}\{e_{(1)}+ie_{(2)}\}.
\]

\begin{Prop}\label{p:nullinv}
The inverse mapping of Proposition \ref{p:adapt} is
\begin{align}
\frac{\partial}{\partial t}&=-|\xi|\sinh r\;e_{(0)}+\frac{|\xi|}{\sqrt{2}}e_{(+)}+\frac{|\xi|}{\sqrt{2}}e_{(-)}, \nonumber\\
\frac{\partial}{\partial z}&=\frac{\bar{\xi}}{2}e_{(0)}-\frac{\bar{\xi}e^{-r}}{2\sqrt{2}}e_{(+)}+\frac{\bar{\xi}e^{r}}{2\sqrt{2}}e_{(-)}.\nonumber
\end{align}
\end{Prop}

\begin{pf}
From Proposition \ref{p:adapt} we have the linear system 
\[
\begin{bmatrix} e_{(0)}\\  \\ e_{(+)}\\  \\ e_{(-)} \end{bmatrix}=A\begin{bmatrix} \frac{\partial}{\partial z}\\  \\ \frac{\partial}{\partial \bar{z}} \\  \\ \frac{\partial}{\partial t} \end{bmatrix},
\] 
where 
\[
A=\frac{1}{\sqrt{2}\xi\bar{\xi}\cosh^2 r}\begin{bmatrix} \xi\sqrt{2} & \bar{\xi}\sqrt{2} & -\sqrt{2}|\xi|\sinh r \\  &   &   \\ -\xi e^{-r} & \bar{\xi}e^{r} & |\xi| \\  &   &   \\ \xi e^{r} & -\bar{\xi}e^{-r} & |\xi| \end{bmatrix}.
\]
So the inverse of $A$ is
\[
A^{-1}=\frac{1}{2\sqrt{2}}\begin{bmatrix} \bar{\xi}\sqrt{2} & -\bar{\xi}e^{-r} & \bar{\xi}e^{r} \\   &    &    \\ \xi\sqrt{2} & \xi e^{r} & -\xi e^{-r} \\  &   &   \\ -2\sqrt{2}|\xi|\sinh r & 2|\xi| & 2|\xi| \end{bmatrix}.
\]
The result follows.
\end{pf}

\begin{Def}
Let $\gamma:{\mathbb{R}}\rightarrow {\mathbb{H}}^3$ be an oriented geodesic. A {\it Jacobi field} along $\gamma$ is a vector field $X$ along $\gamma$ satisfying the equation
\[
\nabla_{\dot{\gamma}}\nabla_{\dot{\gamma}}X^j-X^j=0.
\]
The solutions of this equation form a 6-dimensional vector space for each oriented geodesic $\gamma$, which we denote by 
${\cal{J}}(\gamma)$. Let ${\cal{J}}^\perp(\gamma)$ be the 4-dimensional vector space of Jacobi fields that are orthogonal to $\gamma$.
\end{Def}

\begin{Def}
For an oriented geodesic $\gamma$ in ${\mathbb{H}}^3$, let $Pr_\gamma:T_\gamma{\mathbb{H}}^3\rightarrow T_\gamma{\mathbb{H}}^3$ be projection onto the plane orthogonal to the geodesic.

Let $h:T_\gamma{\mathbb{U}}\rightarrow T_\gamma{\mathbb{H}}^3$ be defined by $h=Pr_\gamma\circ D\Phi$.

\end{Def}

\begin{Prop}
The map $h$ has local coordinate description on ${\mathbb{U}}$
\[
h_{(\xi,\eta,r)}\left(\frac{\partial}{\partial \xi}\right)=-\frac{e^{r}}{2\sqrt{2}\xi}e_{(+)}-\frac{e^{-r}}{2\sqrt{2}\xi}e_{(-)},
\quad
h_{(\xi,\eta,r)}\left(\frac{\partial}{\partial \eta}\right)=-\frac{\bar{\xi}e^{-r}}{2\sqrt{2}}e_{(+)}+\frac{\bar{\xi}e^{r}}{2\sqrt{2}}e_{(-)}.
\]
\end{Prop}
\begin{pf}
From Propositions \ref{p:dphi} and \ref{p:nullinv} we get 

\begin{align}
D\Phi_{(\xi,\eta,r)}\left(\frac{\partial}{\partial \xi}\right)&=-\frac{e^{r}}{2\sqrt{2}\xi}e_{(+)}-\frac{e^{-r}}{2\sqrt{2}\xi}e_{(-)}, \nonumber \\
D\Phi_{(\xi,\eta,r)}\left(\frac{\partial}{\partial \eta}\right)&=\frac{\bar{\xi}}{2}e_{(0)}-\frac{\bar{\xi}e^{-r}}{2\sqrt{2}}e_{(+)}+\frac{\bar{\xi}e^{r}}{2\sqrt{2}}e_{(-)}, \nonumber
\end{align}
and therefore the projection of $D\Phi_{(\xi,\eta,r)}\left(\frac{\partial}{\partial \xi}\right)$ and $D\Phi_{(\xi,\eta,r)}\left(\frac{\partial}{\partial \eta}\right)$ gives the expressions of the Proposition. 
\end{pf}

The next Proposition shows that the tangent $T_{(\xi,\eta)}{\mathbb{U}}$  can be identified with the
orthogonal Jacobi fields along the geodesic ($\xi,\eta$).

\begin{Prop}
The map $h$ is a vector space isomorphism between $T_\gamma{\mathbb{U}}$ and ${\cal{J}}^\perp(\gamma)$.
\end{Prop}
\begin{pf}
First we show that $e_{(+)}$ is a parallel vector field along the geodesic, that is 
\[
\nabla_{e_{(0)}}e_{(+)=0}.
\]
To do this we compute the non-vanishing covariant derivatives:
\begin{align}
\nabla_{\frac{\partial}{\partial t}}\frac{\partial}{\partial t}&=-\frac{1}{2}|\xi|\cosh r\frac{\partial}{\partial t} ,  \qquad \nabla_{\frac{\partial}{\partial z}}\frac{\partial}{\partial \bar{z}}=\frac{1}{2}|\xi|\cosh r\frac{\partial}{\partial t}, \nonumber \\
\nabla_{\frac{\partial}{\partial t}}\frac{\partial}{\partial z}&=-\frac{1}{2}|\xi|\cosh r\frac{\partial}{\partial z}, \qquad \nabla_{\frac{\partial}{\partial t}}\frac{\partial}{\partial \bar{z}}=-\frac{1}{2}|\xi|\cosh r\frac{\partial}{\partial \bar{z}}. \nonumber 
\end{align}
Thus
\[
\nabla_{e_{(0)}}\frac{\partial}{\partial z}=\tanh r \frac{\partial}{\partial z}+\frac{|\xi|}{2\xi\cosh r}\frac{\partial}{\partial t}, \qquad \nabla_{e_{(0)}}\frac{\partial}{\partial \bar{z}}=\tanh r \frac{\partial}{\partial \bar{z}}+\frac{|\xi|}{2\bar{\xi}\cosh r}\frac{\partial}{\partial t} ,
\]
\[
\nabla_{e_{(0)}}\frac{\partial}{\partial t}=-\frac{|\xi|}{\bar{\xi}\cosh r}\frac{\partial}{\partial z}-\frac{|\xi|}{\xi\cosh r}\frac{\partial}{\partial \bar{z}}+\tanh r\frac{\partial}{\partial t} ,
\]
and finally we find that
\[
\nabla_{e_{(0)}}e_{(+)}=0.
\]
Now, since the frames are parallel, for $\alpha=-,+$ and any function $f(r)$
\begin{align}
\nabla_{e_{(0)}}\left(fe_{(\alpha)}\right)&=\frac{df}{dr}e_{(\alpha)}+f\nabla_{e_{(0)}}e_{(\alpha)}=\frac{df}{dr}e_{(\alpha)}, \nonumber \\
\nabla_{e_{(0)}}^2 \left(fe_{(\alpha)}\right)&=\nabla_{e_{(0)}}\left(\frac{df}{dr}e_{(\alpha)}\right)=\frac{d^2 f}{dr^2}e_{(\alpha)}+\frac{df}{dr}\nabla_{e_{(0)}}e_{(\alpha)}=\frac{d^2 f}{dr^2}e_{(\alpha)}. \nonumber
\end{align}
Now we are ready to prove our claim. Since the frame is adapted, we have that $\dot{\gamma}=e_{(0)}$ and
\begin{align}
\nabla_{e_{(0)}}^2 h\left(\frac{\partial}{\partial \xi}\right)&=\nabla_{e_{(0)}}^2\left(-\frac{e^{r}}{2\sqrt{2}\xi}e_{(+)}-\frac{e^{-r}}{2\sqrt{2}\xi}e_{(-)}\right) ,\nonumber \\
&=\frac{d^2}{dr^2}\left(-\frac{e^{r}}{2\sqrt{2}\xi}\right)e_{(+)}+\frac{d^2}{dr^2}\left(-\frac{e^{-r}}{2\sqrt{2}\xi}\right)e_{(-)} ,\nonumber \\
&=-\frac{e^{r}}{2\sqrt{2}\xi}e_{(+)}-\frac{e^{-r}}{2\sqrt{2}\xi}e_{(-)}, \nonumber \\
&=h\left(\frac{\partial}{\partial \xi}\right).\nonumber
\end{align}
Thus $h(\frac{\partial}{\partial \xi})$ is a Jacobi field along the geodesic. Similarly 
$h(\frac{\partial}{\partial \eta})$ can be shown to be a Jacobi field along the geodesic. Moreover, 
these vector fields span the space of orthogonal Jacobi fields along the geodesic.
\end{pf}

\vspace{0.2in}

\section{ The K\"ahler Structure on ${\mathbb{L}}({\mathbb{H}}^3)$}

A K\"ahler structure on a 4-manifold ${\mathbb{M}}$ is a triple (${\mathbb{J}},\Omega,{\mathbb{G}}$), where
${\mathbb{J}}$ is a complex structure, $\Omega$ is a symplectic 2-form and ${\mathbb{G}}$ is an inner product. These
are required to satisfy the conditions:
\[
\Omega({\mathbb{J}}\cdot,{\mathbb{J}}\cdot)=\Omega(\cdot,\cdot),
\qquad\qquad {\mathbb{G}}(\cdot,\cdot)=\Omega({\mathbb{J}}\cdot,\cdot).
\]

We now construct a K\"ahler structure on ${\mathbb{L}}({\mathbb{H}}^3)$. We first define the structure on the open
subset ${\mathbb{U}}\subset{\mathbb{L}}({\mathbb{H}}^3)$ and then show that it extends to the whole space.

\begin{Def}
Given an oriented geodesic $\gamma$ in ${\mathbb{H}}^3$, let 
${\cal{R}}_\gamma:T_\gamma{\mathbb{H}}^3\rightarrow T_\gamma{\mathbb{H}}^3$ be rotation through 90$^{0}$ about the
tangent vector to the geodesic.
\end{Def}

Because we are working in a space of constant curvature this rotation preserves Jacobi fields:

\begin{Prop}\cite{hitch}
The map ${\cal{R}}_\gamma$ takes ${\cal{J}}^\perp$ to ${\cal{J}}^\perp$.
\end{Prop}

We now define our complex structure:

\begin{Def}
Let ${\mathbb{J}}:T_\gamma{\mathbb{U}}\rightarrow T_\gamma{\mathbb{U}}$ be defined to be 
${\mathbb{J}}=h^{-1}\circ{\cal{R}}_\gamma\circ h$.
\end{Def}

It is clear that ${\mathbb{J}}^2=-$Id, so that we have an almost complex structure. In order to be a complex
structure ${\mathbb{J}}$ must also satisfy a certain integrability condition. This is equivalent to the existence of
holomorphic coordinates, which we demonstrate below.

\begin{Prop}
The following two vectors form an eigenbasis for ${\mathbb{J}}$ at ($\xi,\eta$):
\[
\bar{\xi}^2\frac{\partial}{\partial \bar{\xi}}+\frac{\partial}{\partial \eta},
\qquad\qquad
\xi^2\frac{\partial}{\partial \xi}-\frac{\partial}{\partial \bar{\eta}}.
\]
\end{Prop}
\begin{pf}
Note that ${\cal{R}}_\gamma(e_{(+)})=ie_{(+)}$ and thus we have 
\[
{\mathbb{J}}\left(\bar{\xi}^2\frac{\partial}{\partial \bar{\xi}}+\frac{\partial}{\partial \eta}\right)=
h^{-1}\circ{\cal{R}}_\gamma\left(-\frac{\bar{\xi}e^{-r}}{\sqrt{2}}e_{(+)}\right)=i\left(\bar{\xi}^2\frac{\partial}{\partial \bar{\xi}}+\frac{\partial}{\partial \eta}\right).
\]
Similarly, we have
\[
{\mathbb{J}}\left(\xi^2\frac{\partial}{\partial \xi}-\frac{\partial}{\partial \bar{\eta}}\right)=i\left(\xi^2\frac{\partial}{\partial \xi}-\frac{\partial}{\partial \bar{\eta}}\right),
\]
as claimed.
\end{pf}

We introduce holomorphic coordinates on ${\mathbb{U}}$ as follows:

\begin{Def}\label{d:mu12}
Define  $(\mu_{1},\mu_{2})\in{\mathbb{C}}^2$ by
\[
\mu_{1}=-\eta+\frac{1}{\bar{\xi}}\; , \qquad \mu_{2}=\frac{1}{\bar{\eta}+\frac{1}{\xi}},
\]
the inverse relation being
\[
\eta=-\frac{1}{2}\left(\mu_{1}-\frac{1}{\bar{\mu}_{2}}\right), \qquad \xi=\frac{2}{\bar{\mu}_{1}+\frac{1}{\mu_{2}}}.
\]
\end{Def}

We are now in a position to prove one of our main results:

\begin{Thm}
The almost complex structure ${\mathbb{J}}$ extends to the whole of ${\mathbb{L}}({\mathbb{H}}^3)$ and is integrable.
Moreover, the complex surface (${\mathbb{L}}({\mathbb{H}}^3),{\mathbb{J}})$ is biholomorphic to 
${\mathbb{P}}^1\times{\mathbb{P}}^1-\bar{\Delta}$.
\end{Thm}
\begin{pf}
We first note that $(\mu_{1},\mu_{2})$ are holomorphic coordinates on ${\mathbb{U}}$. Indeed,
\begin{align}
\frac{\partial}{\partial \mu_{1}}&=-\frac{1}{2}\left[\frac{2\bar{\mu}_{2}}{(\mu_{1}\bar{\mu}_{2}+1)}\right]^2\frac{\partial}{\partial \bar{\xi}}-\frac{1}{2}\frac{\partial}{\partial \eta}=-\frac{1}{2}\left(\bar{\xi}^2 \frac{\partial}{\partial \bar{\xi}}+\frac{\partial}{\partial \eta}\right), \nonumber \\
\frac{\partial}{\partial \mu_{2}}&=\frac{1}{2\mu_{2}^2}\left[\left( \frac{2\mu_{2}}{\bar{\mu}_{1}\mu_{2}+1}\right)^2 \frac{\partial}{\partial \xi}-\frac{\partial}{\partial \bar{\eta}}\right]=\frac{1}{2}\left(\bar{\eta}+\frac{1}{\xi}\right)^2 \left(\xi^2 \frac{\partial}{\partial \xi}-\frac{\partial}{\partial \bar{\eta}}\right), \nonumber 
\end{align}
and hence we find that
\[
{\mathbb{J}}\left(\frac{\partial}{\partial \mu_{1}}\right)=i\frac{\partial}{\partial \mu_{1}}, \qquad {\mathbb{J}}\left(\frac{\partial}{\partial \mu_{2}}\right)=i\frac{\partial}{\partial \mu_{2}},
\]
so we have holomorphic coordinates for ${\mathbb{J}}$ and therefore the almost complex structure is integrable. 

These coordinates come from the past and future boundaries of the oriented geodesic in the ball model. To see this 
explicitly, introduce coordinates $(w,\rho)$ on the 3-ball by $w=y^{1}+iy^{2}$ and $\rho=y^{3}$. The isometry $P$ 
given by equations (\ref{e:ballcoords}) becomes 
$P(z,t)=(w,\rho)$, with
\[
w=\frac{2z}{z\bar{z}+(t+1)^2}\; , \qquad \rho=1-\frac{2(t+1)}{z\bar{z}+(t+1)^2}.
\]
Mapping a geodesic $\gamma\in{\mathbb{U}}$ to the Poincar\'e ball we get
\[
w=\frac{2\left[\left(\eta+\frac{1}{\bar{\xi}}\right)e^{2r}+ \eta-\frac{1}{\bar{\xi}}\right]}{\left(\eta\bar{\eta}+\frac{1}{\xi\bar{\xi}}+\frac{\eta}{\xi}+\frac{\bar{\eta}}{\bar{\xi}}+1\right)e^{2r}+\frac{4}{|\xi|}e^{r}+\left(\eta\bar{\eta}+\frac{1}{\xi\bar{\xi}}-\frac{\eta}{\xi}-\frac{\bar{\eta}}{\bar{\xi}}+1    \right)},
\]
\[
\rho=1-\frac{2\left(e^{2r}+\frac{2}{|\xi|}e^{r}+1\right)}{\left(\eta\bar{\eta}+\frac{1}{\xi\bar{\xi}}+\frac{\eta}{\xi}+\frac{\bar{\eta}}{\bar{\xi}}+1\right)e^{2r}+\frac{4}{|\xi|}e^{r}+\left(\eta\bar{\eta}+\frac{1}{\xi\bar{\xi}}-\frac{\eta}{\xi}-\frac{\bar{\eta}}{\bar{\xi}}+1    \right)}.
\]

Define 
\[
w_{\pm}=\lim_{r\to\pm\infty}w(r) \qquad \rho_{\pm}=\lim_{r\to\pm\infty}\rho(r),
\]
and find 
\[
w_{+}=\frac{2\left(\eta+\frac{1}{\bar{\xi}}\right)}{\eta\bar{\eta}+\frac{1}{\xi\bar{\xi}}+\frac{\eta}{\xi}+\frac{\bar{\eta}}{\bar{\xi}}+1}=\frac{2\mu_{2}}{\mu_{2}\bar{\mu}_{2}+1},
\]
\[
\rho_{+}=1-\frac{2}{\eta\bar{\eta}+\frac{1}{\xi\bar{\xi}}+\frac{\eta}{\xi}+\frac{\bar{\eta}}{\bar{\xi}}+1}=\frac{1-\mu_{2}\bar{\mu}_{2}}{1+\mu_{2}\bar{\mu}_{2}},
\]
\[
w_{-}=\frac{2\left(\eta-\frac{1}{\bar{\xi}}\right)}{\eta\bar{\eta}+\frac{1}{\xi\bar{\xi}}-\frac{\eta}{\xi}-\frac{\bar{\eta}}{\bar{\xi}}+1}=-\frac{2\mu_{1}}{\mu_{1}\bar{\mu}_{1}+1},
\]
\[
\rho_{-}=1-\frac{2}{\eta\bar{\eta}+\frac{1}{\xi\bar{\xi}}-\frac{\eta}{\xi}-\frac{\bar{\eta}}{\bar{\xi}}+1}=\frac{\mu_{1}\bar{\mu}_{1}-1}{\mu_{1}\bar{\mu}_{1}+1}.
\]

We see that $w_{\pm}\overline{w}_{\pm}+\rho_{\pm}^2=1$ and so $(w_{\pm},\rho_{\pm})\in \partial B^3=S^2$, as expected.
In fact, $\mu_2$ is the holomorphic coordinate obtained on S$^2$ by stereographic projection from the south pole,
while $\mu_1$ is the anti-holomorphic coordinate (the composition of stereographic projection and the antipodal
map $\iota$).

It is now clear that the definition of ${\Bbb{J}}$ as rotation of Jacobi fields about the oriented geodesic 
extends to all of ${\mathbb{L}}({\mathbb{H}}^3)$ and the Theorem follows.  
\end{pf}

We now define the symplectic structure on ${\mathbb{U}}$.

\begin{Def}
Let ${\mathbb{X}},{\mathbb{Y}}\in T_\gamma{\mathbb{U}}$. Define a 2-form on ${\mathbb{U}}$ by
\[
\Omega({\mathbb{X}},{\mathbb{Y}})=g(h({\mathbb{X}}),\nabla_{\dot{\gamma}}(h({\mathbb{Y}}))
    -g(h({\mathbb{Y}}),\nabla_{\dot{\gamma}}(h({\mathbb{X}})).
\]
\end{Def}

\begin{Prop}\label{p:symp1}
The 2-form $\Omega$ is given in local coordinates $(\xi,\eta)$ by
\[
\Omega=-\frac{1}{2}\{ d\xi\wedge d\bar{\eta}+d\bar{\xi}\wedge d\eta\},
\]
and so is closed and non-degenerate - that is, it is a symplectic structure.  
\end{Prop}
\begin{pf}
We have found that
\begin{align}
h\left(\frac{\partial}{\partial \xi}\right)=-\frac{e^{r}}{2\sqrt{2}\xi}e_{(+)}-\frac{e^{-r}}{2\sqrt{2}\xi}e_{(-)}, \quad \nabla_{\dot{\gamma}}h\left(\frac{\partial}{\partial \xi}\right)&=-\frac{e^{r}}{2\sqrt{2}\xi}e_{(+)}+\frac{e^{-r}}{2\sqrt{2}\xi}e_{(-)}, \nonumber \\
h\left(\frac{\partial}{\partial \eta}\right)=-\frac{\bar{\xi}e^{-r}}{2\sqrt{2}}e_{(+)}+\frac{\bar{\xi}e^{r}}{2\sqrt{2}}e_{(-)}, \quad \nabla_{\dot{\gamma}}h\left(\frac{\partial}{\partial \eta}\right)&=\frac{\bar{\xi}e^{-r}}{2\sqrt{2}}e_{(+)}+\frac{\bar{\xi}e^{r}}{2\sqrt{2}}e_{(-)}. \nonumber
\end{align}
Thus, the only non-vanishing components of $\Omega$ can be computed to be
\[
\Omega\left(\frac{\partial}{\partial \xi},\frac{\partial}{\partial \bar{\eta}}\right)=\Omega\left(\frac{\partial}{\partial \bar{\xi}},\frac{\partial}{\partial \eta}\right)=-\frac{1}{2}.
\]
\end{pf}

The complex structure and symplectic form are compatible in the following sense:

\begin{Prop}
\[
\Omega({\mathbb{J}}\cdot,{\mathbb{J}}\cdot)=\Omega(\cdot,\cdot).
\]
\end{Prop}

\begin{pf}
Let ${\mathbb{X}},{\mathbb{Y}}\in T_\gamma{\mathbb{U}}$ and we compute
\begin{align}
\Omega({\mathbb{J}}({\mathbb{X}}),{\mathbb{J}}({\mathbb{Y}}))=&g(h({\mathbb{J}}({\mathbb{X}})),\nabla_{\dot{\gamma}}(h({\mathbb{J}}({\mathbb{Y}})))
    -g(h({\mathbb{J}}({\mathbb{Y}})),\nabla_{\dot{\gamma}}(h({\mathbb{J}}({\mathbb{X}}))),\nonumber\\
&=g({\cal{R}}_\gamma(h({\mathbb{X}})),\nabla_{\dot{\gamma}}({\cal{R}}_\gamma(h({\mathbb{Y}})))
    -g({\cal{R}}_\gamma(h({\mathbb{Y}})),\nabla_{\dot{\gamma}}({\cal{R}}_\gamma(h({\mathbb{X}})))).\nonumber
\end{align}
It is clear that the Proposition holds if $\nabla_{\dot{\gamma}}$ and ${\cal{R}}_\gamma$ commute on orthogonal vector 
fields along the geodesic. To see this,
let $X=f\;e_{(+)}+\bar{f}\;e_{(-)}$ be an orthogonal vector field along the geodesic $\gamma$. Then
\begin{align}
\nabla_{\dot{\gamma}}({\cal{R}}_\gamma(X))=&\nabla_{\dot{\gamma}}i(f\;e_{(+)}-\bar{f}\;e_{(-)})
   =i(\dot{f}\;e_{(+)}-\dot{\bar{f}}\;e_{(-)}),\nonumber\\
&={\cal{R}}_\gamma(\dot{f}\;e_{(+)}+\dot{\bar{f}}\;e_{(-)})
   ={\cal{R}}_\gamma(\nabla_{\dot{\gamma}}(X)),\nonumber
\end{align}
as claimed.
\end{pf}

\begin{Def}
The K\"ahler metric ${\mathbb{G}}$ on ${\mathbb{U}}$ is defined by
\[
{\mathbb{G}}(\cdot\; ,\cdot)=\Omega({\mathbb{J}}\cdot\; ,\cdot).
\]
\end{Def}

This has local coordinate expression:

\begin{Prop}\label{p:kaehlermetric}
The K\"ahler metric has the following expression in local coordinates $(\xi,\eta)$:
\[
{\mathbb{G}}=-\frac{i}{4}\left(\frac{1}{\xi^2}d\xi^2-\frac{1}{\bar{\xi}^2}d\bar{\xi}^2+\bar{\xi}^2d\eta^2-\xi^2d\bar{\eta}^2\right).
\]
\end{Prop}
\begin{pf}
First we express the symplectic form $\Omega$ in terms of $\mu_{1},\mu_{2}$. By Definition \ref{d:mu12} we have
\[
d\eta =-\frac{1}{2}\left(d\mu_{1}+\frac{1}{\bar{\mu}_2^2}d\bar{\mu}_{2}\right), \qquad  d\xi=-\frac{2}{(1+\bar{\mu}_{1}\mu_{2})^2}(\mu_{2}^2 d\bar{\mu}_{1}-d\mu_{2}).
\]
By Proposition \ref{p:symp1} we have $\Omega=-{\mathbb{R}}\mbox{e}[d\xi\wedge d\bar{\eta}]$ and since
\begin{align}
d\xi\wedge d\bar{\eta}&=\frac{1}{(1+\bar{\mu}_{1}\mu_{2})^2}\left(\mu_{2}^2 d\bar{\mu}_{1}-d\mu_{2}\right)\wedge \left(d\bar{\mu}_{1}+\frac{1}{\mu_{2}^2}d\mu_{2}\right) ,\nonumber \\
&=\frac{2}{(1+\bar{\mu}_{1}\mu_{2})^2}d\bar{\mu}_{1}\wedge d\mu_{2}, \nonumber
\end{align}
we get
\begin{equation}\label{e:symp}
\Omega=-\left[\frac{1}{(1+\mu_{1}\bar{\mu}_{2})^2}d\mu_{1}\wedge d\bar{\mu}_{2}+\frac{1}{(1+\bar{\mu}_{1}\mu_{2})^2}d\bar{\mu}_{1}\wedge d\mu_{2}\right].
\end{equation}
Now we find ${\mathbb{G}}$ in these coordinates
\begin{align}
{\mathbb{G}}&=-i\left[\frac{1}{(1+\mu_{1}\bar{\mu}_{2})^2}d\mu_{1}\otimes d\bar{\mu}_{2}-\frac{1}{(1+\bar{\mu}_{1}\mu_{2})^2}d\bar{\mu}_{1}\otimes d\mu_{2}\right], \label{e:metmu}\\
&=-i\frac{(\eta\bar{\xi}+1)^2}{4}\frac{1}{2}\left(d\eta +\frac{1}{\bar{\xi}^2}d\bar{\xi}\right)\otimes \frac{2}{(\eta\bar{\xi}+1)^2}(\bar{\xi}^2 d\eta-d\bar{\xi}), \nonumber \\
& \qquad \qquad +i\frac{(\bar{\eta}\xi+1)^2}{4}\frac{1}{2}\left(d\bar{\eta} +\frac{1}{\xi^2}d\xi\right)\otimes \frac{2}{(\bar{\eta}\xi+1)^2}(\xi^2 d\bar{\eta}-d\xi), \nonumber \\
&=-i\left[\frac{1}{4}(d\eta+\frac{1}{\bar{\xi}^2}d\bar{\xi})\otimes(\bar{\xi}^2d\eta-d\bar{\xi})-\frac{1}{4}(d\bar{\eta}+\frac{1}{\xi^2}d\xi)\otimes(\xi^2 d\bar{\eta}-d\xi)\right], \nonumber \\
&=-\frac{i}{4}\left(\frac{1}{\xi^2}d\xi^2-\frac{1}{\bar{\xi}^2}d\bar{\xi}^2+\bar{\xi}^2d\eta^2-\xi^2d\bar{\eta}^2\right), \nonumber
\end{align}
as claimed. 
\end{pf}

We are now in a position to prove our second result:

\begin{Thm}
The K\"ahler structure is defined on the whole of the space ${\mathbb{L}}({\mathbb{H}}^3)$. The metric ${\mathbb{G}}$
is of neutral signature, is conformally flat and scalar flat.
\end{Thm}
\begin{pf}
It is clear from the expression of ${\mathbb{G}}$ and $\Omega$ in holomorphic coordinates (equations (\ref{e:symp}) and (\ref{e:metmu}))
that these are well defined except where $\mu_1=-\mu_2^{-1}$. But this is just the reflected diagonal, and so the
K\"ahler structure is well-defined on the whole of the space ${\mathbb{L}}({\mathbb{H}}^3)$.

The signature of the metric is ($++--$) and the curvature can be computed directly from the coordinate expressions given 
above. The only non-vanishing components of the Riemann curvature tensor are
\[
R_{\bar{\mu}_1\mu_2\mu_2}^{\;\;\;\;\;\;\;\;\;\;\;\;\mu_2}
=-R_{\bar{\mu}_1\mu_2\bar{\mu}_1}^{\;\;\;\;\;\;\;\;\;\;\;\;\bar{\mu}_1}=\frac{2}{(1+\bar{\mu}_1\mu_2)^2}.
\]
The results are then as stated.
\end{pf}

As in the case of oriented lines of ${\mathbb{E}}^3$, this metric has the following mechanical interpretation:
the length of a vector ${\mathbb{X}}\in T_\gamma{\mathbb{L}}({\mathbb{H}}^3)$ is the angular momentum of the
Jacobi field $h({\mathbb{X}})$ about the geodesic $\gamma$ in ${\mathbb{H}}^3$.

\vspace{0.2in}

\section{The Isometry Group of the Neutral K\"ahler Metric}

We now find the isometry group of the space of oriented geodesics endowed with the above K\"ahler metric.

\begin{Thm}
The identity component of the isometry group of the metric ${\mathbb{G}}$ on ${\mathbb{L}}({\mathbb{H}}^3)$ 
is isomorphic to the identity component of the hyperbolic isometry group.  
\end{Thm}

We prove this by way of a number of propositions: first we find the Killing vectors of ${\mathbb{G}}$
and of the hyperbolic metric $g$. Then we integrate those of $g$ up to get explicit expressions for the isometry
group of ${\mathbb{H}}^3$ and find how this group acts on ${\mathbb{L}}({\mathbb{H}}^3)$. Finally we show that
the infinitesimal generators of this action are precisely the Killing vectors of ${\mathbb{G}}$.

We start then with:

\begin{Prop}\label{p:lieLH}
The Killing vectors of ${\mathbb{G}}$ form a 6-parameter Lie algebra given by
\[
{\mathbb{K}}={\mathbb{R}}{\mbox e}\left[(c_1 +c_2 \mu_{1}+c_3\mu_{1}^2)\frac{\partial}{\partial \mu_{1}}
+ (\bar{c}_3 -\bar{c}_2\mu_{2}+\bar{c}_1\;\mu_{2}^2)\frac{\partial}{\partial \mu_{2}}\right],
\]
where $c_1,c_2,c_3\in {\mathbb{C}}$.
\end{Prop}
\begin{pf}
Let
\[
{\mathbb{K}}={\mathbb{R}}{\mbox e}\left({\mathbb{K}}^{\mu_1}\frac{\partial}{\partial \mu_{1}}
+ {\mathbb{K}}^{\mu_2}\frac{\partial}{\partial \mu_{2}}\right),
\]
be a vector field on ${\mathbb{L}}({\mathbb{H}}^3)$ with  
${\mathbb{K}}^{\mu_i}={\mathbb{K}}^{\mu_i}(\mu_{1},\bar{\mu}_{1},\mu_{2},\bar{\mu}_{2})$ for $i=1,2$.

We will solve the Killing equations for ${\mathbb{K}}$:
\begin{equation}\label{e:killLH}
{\mathbb{K}}^{i}\partial_{i}{\mathbb{G}}_{jk}+{\mathbb{G}}_{ki}\partial_{j}{\mathbb{K}}^{i}+{\mathbb{G}}_{ji}\partial_{k}{\mathbb{K}}^{i}=0,
\end{equation}
where $\partial$ represents partial differentiation.

In what follows we denote the Killing equation with specific values of the indices $j$ and $k$ simply by $(j,k)$. Thus, 
for example, 
$(\mu_{1},\bar{\mu}_{1})$ will mean the Killing equation (\ref{e:killLH}) with $j=\mu_{1},  k=\bar{\mu}_{1}$.
Similarly for the derivative of these equations.

To start then, the $(\bar{\mu}_{1},\bar{\mu}_{1})$ and $(\bar{\mu}_{2},\bar{\mu}_{2})$ equations are
\[
\partial_{\bar{\mu}_{1}}{\mathbb{K}}^{\mu_2}=0, \qquad \partial_{\bar{\mu}_{2}}{\mathbb{K}}^{\mu_1}=0,
\]
and so
\[
{\mathbb{K}}^{\mu_1}={\mathbb{K}}^{\mu_1}(\mu_{1},\bar{\mu}_{1},\mu_{2}), 
\qquad {\mathbb{K}}^{\mu_2}={\mathbb{K}}^{\mu_2}(\mu_{1},\mu_{2},\bar{\mu}_{2}).
\]
Now, differentiating the $(\mu_{1},\mu_{2})$ equation with respect to $\mu_{2}$ gives
\[
(\mu_{2}\bar{\mu}_{1}+1)^2\partial_{\mu_{2}}^2 \bar{{\mathbb{K}}}^{\mu_2}+2(\mu_{2}\bar{\mu}_{1}+1)\partial_{\mu_{2}}\bar{{\mathbb{K}}}^{\mu_2}=0.
\]
This can be rewritten
\[
\partial_{\mu_{2}}\left[(\mu_{2}\bar{\mu}_{1}+1)^2\partial_{\mu_{2}}\bar{{\mathbb{K}}}^{\mu_2}\right]=0,
\] 
which integrates to 
\[
{\mathbb{K}}^{\mu_2}=\alpha_{2}+\frac{\alpha_{1}}{\mu_{1}\bar{\mu}_{2}+1}, \qquad \mbox{ for}\qquad \alpha_{i}=\alpha_{i}(\mu_{1},\bar{\mu}_{1},\mu_{2})\in {\mathbb{C}}, \quad i=1,2.
\]
Since $\partial_{\bar{\mu}_{1}}{\mathbb{K}}^{\mu_2}=0$ we get
\[
(\mu_{1}\bar{\mu}_{2}+1)\partial_{\bar{\mu}_{1}}\alpha_{2}+\partial_{\bar{\mu}_{1}}\alpha_{1}=0.
\]
Differentiating this with respect to $\bar{\mu}_{2}$ we see that $\partial_{\bar{\mu}_{1}}\alpha_{2}=0$, and
hence $\partial_{\bar{\mu}_{1}}\alpha_{1}=0$.
Thus $\alpha_{1}$ and $\alpha_{2}$ are holomorphic, that is $\alpha_{1}=\alpha_{1}(\mu_{1},\mu_{2}), \; \; \alpha_{2}=\alpha_{2}(\mu_{1},\mu_{2})$.

We now differentiate the equation $(\mu_{1},\mu_{2})$ with respect to $\mu_{1}$ yielding 
\[
(\mu_{1}\bar{\mu}_{2}+1)^2\partial_{\mu_{1}}^2 \bar{{\mathbb{K}}}^{\mu_1}+2(\mu_{1}\bar{\mu}_{2}+1)\partial_{\mu_{1}}\bar{{\mathbb{K}}}^{\mu_1}=0,
\]
which, by a similar argument, has solution
\[
{\mathbb{K}}^{\mu_1}=\alpha_{3}+\frac{\alpha_{4}}{\bar{\mu}_{1}\mu_{2}+1}, \qquad \mbox{ for}\qquad \alpha_{3}(\mu_{1},\mu_{2}),\; \;\alpha_{4}(\mu_{1},\mu_{2})\in {\mathbb{C}}.
\]
The equation $(\bar{\mu}_{1},\bar{\mu}_{2})$ now gives $\alpha_{1}=\alpha_{4}=0$, while $\partial_{\mu_{1}} \partial_{\bar{\mu}_{2}}^2 (\mu_{1},\bar{\mu}_{1})$ says that
\[
\mu_{1}\partial_{\mu_{1}}^2 \alpha_{2}+2\partial_{\mu_{1}} \alpha_{2}=0,
\]
with solution
\[
\alpha_{2}=\alpha_{5}+\frac{\alpha_{6}}{\mu_{1}}, \qquad\mbox{ for}\qquad \alpha_{5}(\mu_{2}),\alpha_{6}(\mu_{2})\in {\mathbb{C}}.
\] 
Now, $\partial_{\bar{\mu}_{1}}^2 \partial_{\mu_{2}} (\mu_{1},\bar{\mu}_{1})$ reads
\[
-(\mu_{1}\bar{\mu}_{2}+1)^2\partial_{\mu_{2}} \alpha_{6}+2\bar{\alpha}_{6}\mu_{1}^2\mu_{2}=0,
\]
which implies that $\alpha_{6}=0$. On the other hand the derived equation 
$\partial_{\bar{\mu}_{1}}^2\partial_{\mu_{2}}(\mu_{2},\bar{\mu}_{2})$ implies that
\[
\mu_{2}\partial_{\mu_{2}}^2\alpha_{3}+2\partial_{\mu_{2}}\alpha_{3}=0,
\]
and so
\[
\alpha_{3}=\alpha_{7}+\frac{\alpha_{8}}{\mu_{2}},
\]
where $\alpha_{7}=\alpha_{7}(\mu_{1}),\alpha_{8}=\alpha_{8}(\mu_{1})\in {\mathbb{C}}$. 

Substituting the previous results in equation $(\mu_{2},\bar{\mu}_{2})$ we find that
\[
(\bar{\mu}_{1}\mu_{2}+1)^2\partial_{\mu_{1}}\alpha_{8}-2\bar{\alpha}_{8}\mu_{1}\mu_{2}^2=0,
\]
from which we conclude that $\alpha_{8}=0$.

The derived equation $\partial_{\bar{\mu}_{2}}\partial_{\mu_{1}}^2(\mu_{1},\bar{\mu}_{2})$
implies that $\partial_{\mu_{1}}^3\alpha_{7}=0$ and so
\[
\alpha_{7}=c_{1}+c_{2}\mu_{1}+c_{3}\mu_{1}^2\; , \qquad \mbox{ for}\qquad c_1,c_2,c_3\in {\mathbb{C}}.
\]
Similarly, $\partial_{\bar{\mu}_{1}}\partial_{\mu_{2}}^2(\bar{\mu}_{1},\mu_{2})$ implies that
$\partial_{\mu_{2}}^3\alpha_{5}=0$ and we therefore have
\[
\alpha_{5}=c_{4}+c_{5}\mu_{2}+c_{6}\mu_{2}^2\; , \qquad \mbox{ for}\qquad c_{4},c_{5},c_{6}\in {\mathbb{C}}.
\]

Finally, putting all of the preceding together we obtain 
\[
{\mathbb{K}}^{\mu_1}=c_{1}+c_{2}\mu_{1}+c_{3}\mu_{1}^2
\qquad \qquad {\mathbb{K}}^{\mu_2}=\bar{c}_{3}-\bar{c}_{2}\mu_{2}+\bar{c}_{1}\mu_{2}^2,
\]
as claimed.
\end{pf}

By a similar method we compute the isometry group of the upper half-space model of ${\mathbb{H}}^3$ and its 
Lie algebra.

\begin{Prop}\label{p:killLH}
The Killing vectors of the hyperbolic metric $g$ form a 6-parameter Lie algebra given by
\[
K=K^0 \frac{\partial}{\partial x^0}+K^1 \frac{\partial}{\partial x^1}+K^2 \frac{\partial}{\partial x_2},
\]
where
\[
K^{0}=(a_{4}+a_{5}x^1+a_{6}x^2)x^0,
\]
\[
K^{1}=a_{11}-{\textstyle{\frac{1}{2}}}a_{5}\left((x^0)^2-(x^1)^2+(x^2)^2\right)-a_{10}x^2+a_{4}x^1+a_{6}x^1x^2,
\]
\[
K^{2}=a_{9}-{\textstyle{\frac{1}{2}}}a_{6}\left((x^0)^2+(x^1)^2-(x^2)^2\right)+a_{10}x^1+a_{4}x^2+a_{5}x^1x^2,
\]
for $a_4,a_5,a_6,a_9,a_{10},a_{11}\in{\mathbb{R}}$.
\end{Prop}
\begin{pf}
Let $K=K^0 \frac{\partial}{\partial x^0}+K^1 \frac{\partial}{\partial x^1}+K^2 \frac{\partial}{\partial x_2}$ be a 
Killing vector of the hyperbolic metric $g$.  Then it satisfies
\begin{equation}\label{e:killH}
K^{i}\partial_{i}g_{jk}+g_{ki}\partial_{j}K^{i}+g_{ji}\partial_{k}K^{i}=0,
\end{equation}
where $g_{ij}=(x^0)^{-2}\delta_{ij}$.

As before, we denote the Killing equation with specific values of the indices $j$ and $k$ simply by $(j,k)$. Thus, 
for example, $(x^1,x^2)$ will mean the Killing equation (\ref{e:killH}) with $j=x^1,  k=x^2$. Throughout, all functions
will be real-valued.

To start then, the $(x^0,x^0)$ equation is
\[
x^0\partial_{x^0}K^0-K^0=0,
\]
with solution $K^0=a_{1}x^0$ for $a_{1}=a_{1}(x^1,x^2)$.
Substituting this into the equations $(x^1,x^0)$ and $(x^2,x^0)$ we obtain
\[
x^0\partial_{x^1}a_{1}+\partial_{x^0}K^1=0 \qquad\qquad x^0\partial_{x^2}a_{1}+\partial_{x^0}K^2,
\]
which we integrate to $K^1=a_{2}-{\textstyle{\frac{1}{2}}}(x^0)^2\partial_{x^1}a_{1}$
and $K^2=a_{3}-{\textstyle{\frac{1}{2}}}(x^0)^2\partial_{x^2}a_{1}$ for $a_{2}=a_{2}(x^1,x^2)$ and
$a_{3}=a_{3}(x^1,x^2)$.

Now the derived equations $\partial_{x^0}^2 (x^2,x^2),\; \partial_{x^0}^2 (x^1,x^1),\; \partial_{x^0}^2 (x^2,x^1)$ yield 
\[
\partial_{x^2}^2 a_{1}=0,\qquad\quad \partial_{x^1}^2 a_{1}=0,\qquad\quad \partial_{x^2}\partial_{x^1} a_{1}=0.
\]
The first of these implies $a_{1}=a+bx^2$ for $a=a(x^1)$ and $b=b(x^1)$, while the last equation means that
$\partial_{x^1}b=0$ and so $b$ is constant. The middle equation says that 
$\partial_{x^1}(\partial_{x^1}a_{1})=\partial_{x^1}^2 a=0$ and hence $a=c+dx^1$, where $c,d$ are constants.
Therefore we have $a_{1}=a_{4}+a_{5}x^1+a_{6}x^2$, where $ a_{4},a_{5},a_{6}$ are constants.

Substituting the above results to the equations $(x^1,x^1)$ and $(x^2,x^2)$ we obtain
\[
\partial_{x^1}a_{2}-a_{4}-a_{5}x^1-a_{6}x^2=0\qquad \partial_{x^2}a_{3}-a_{4}-a_{5}x^1-a_{6}x^2=0,
\]
with solutions 
\[
a_{2}=a_{7}+a_{4}x^1+{\textstyle{\frac{1}{2}}}a_{5}(x^1)^2+a_{6}x^1x^2,
\]
and 
\[
a_{3}=a_{8}+a_{4}x^2+{\textstyle{\frac{1}{2}}}a_{6}(x^2)^2+a_{5}x^1x^2,
\]
where $a_{7}=a_{7}(x^2)$ and $a_{8}=a_{8}(x^1)$.

Now the equation $\partial_{x^1}(x^2,x^1)$ gives $\partial_{x^1}^2 a_{8}+a_{6}=0$ and therefore
\[
a_{8}=a_{9}+a_{10}x^1-{\textstyle{\frac{1}{2}}}a_{6}(x^1)^2,
\]
where $a_{9},a_{10}$ are constants.

From $\partial_{x^2}(x^2,x^1)$ says that $\partial_{x^2}^2 a_{7}+a_{5}=0$ and therefore
\[
a_{7}=a_{11}+a_{12}x^2-{\textstyle{\frac{1}{2}}}a_{5}(x^2)^2,
\]
where $a_{11},a_{12}$ are constants.
The equation $(x^2,x^1)$ gives $a_{10}=-a_{12}$.

Finally, assembling the expressions for $K^{i}$ we get:
\begin{align}
K^{0}&=(a_{4}+a_{5}x^1+a_{6}x^2)x^0,\nonumber \\
K^{1}&=a_{11}-{\textstyle{\frac{1}{2}}}a_{5}\left((x^0)^2-(x^1)^2+(x^2)^2\right)-a_{10}x^2+a_{4}x^1+a_{6}x^1x^2,\nonumber \\
K^{2}&=a_{9}-{\textstyle{\frac{1}{2}}}a_{6}\left((x^0)^2+(x^1)^2-(x^2)^2\right)+a_{10}x^1+a_{4}x^2+a_{5}x^1x^2,\nonumber 
\end{align}
as claimed.
\end{pf}
 
\begin{Cor}
The Killing vectors of the hyperbolic metric $g$ are
\[
K={\mathbb{R}}\mbox{e}\left[t(\gamma+2\alpha\bar{z})\frac{\partial}{\partial t}
   +2\left(\beta+\gamma z-\alpha t^2+\bar{\alpha}z^2\right)\frac{\partial}{\partial z}\right],
\]
for $\alpha,\beta,\gamma\in{\mathbb{C}}$.
\end{Cor}
\begin{pf}
If we re-introduce variables $(t,z)$ in the upper half-space model of ${\mathbb{H}}^3$ by
\[
t=x^0 \in {\mathbb{R}}_+, \quad z=x^1+ix^2\in {\mathbb{C}},
\]
and set
\[
\alpha={\textstyle{\frac{1}{2}}}(a_{5}+ia_{6}), \quad \beta=a_{11}+ia_{9}, \quad \gamma=a_{4}+ia_{10}, 
\qquad \alpha,\beta,\gamma\in {\mathbb{C}},
\]
the result follows from Proposition \ref{p:killLH}.
\end{pf}

We now integrate the Killing vectors of ${\mathbb{H}}^3$ to get the group action:

\begin{Prop}
The isometry group of ${\mathbb{H}}^3$ is 6-dimensional and the associated 1-parameter group of local isometries which map
$((t_0,z_0),s)\rightarrow(t(s),z(s))$, for $\alpha$ and  $\gamma$ not both zero, are
\[
t(s)=\frac{t_0(|z_{0}-\tau|^2+t_{0}^2)}
  {\left|\bar{z}_0-\bar{\tau}-\frac{\bar{\alpha}}{\gamma_{1}}( e^{\gamma_{1}s}-1)(|z_{0}-\tau|^2+t_{0}^2)\right|^2+t_0^2}
  e^{\frac{\gamma_{1}+\bar{\gamma}_{1}}{2}s},
\]
\[
z(s)=\frac{z_0-\tau-\frac{\alpha}{\bar{\gamma}_{1}}(e^{\bar{\gamma}_{1}s}-1)(|z_{0}-\tau|^2+t_{0}^2)}
 {\left|\bar{z}_0-\bar{\tau}-\frac{\bar{\alpha}}{\gamma_{1}}( e^{\gamma_{1}s}-1)(|z_{0}-\tau|^2+t_{0}^2)\right|^2
     +t_0^2}(|z_{0}-\tau|^2+t_{0}^2)e^{\gamma_{1}s}+\tau ,
\]
where $\beta=-\gamma\tau-\bar{\alpha}\tau^2$ and $\gamma_1=\gamma+2\bar{\alpha}\tau$.

For $\alpha=\gamma=0$ the isometries are
\[
t(s)=t_{0},\qquad z(s)=\beta s+z_{0}.
\]
\end{Prop}
\begin{pf}
Consider the integral curve ${\mathbb{R}}\rightarrow {\mathbb{H}}^3\colon s \mapsto t(s)\frac{\partial}{\partial t}+z(s)\frac{\partial}{\partial z}$, of the vector field $K$. To find the integral curves explicitly of $K$ we have to 
solve the system of differential equations
\begin{equation}\label{e:lieH}
t'=t\;\left({\textstyle{\frac{1}{2}}}(\gamma+\bar{\gamma})+\alpha\bar{z}+\bar{\alpha}z\right), \qquad z'=\beta+\gamma z-\alpha t^2+\bar{\alpha}z^2,
\end{equation}
where a prime denotes differentiation with respect of the variable $s$.

For $\alpha=\gamma=0$ the result follows immediately, so we now consider the case where $\alpha$ and  $\gamma$ not both 
zero.

First let us assume that $\beta=0$. The system (\ref{e:lieH}) can be written as
\begin{equation}\label{e:mateq}
V'=VAV+BV+VB^{T},
\end{equation}
where
\[
V=\begin{bmatrix} z & t \\ -t & \bar{z} \end{bmatrix}, \quad A=\begin{bmatrix} \bar{\alpha} & 0 \\ 0 & \alpha \end{bmatrix}, \quad B={\textstyle{\frac{1}{2}}}\begin{bmatrix} \gamma & 0 \\ 0 & \bar{\gamma} \end{bmatrix} .
\]
The system (\ref{e:mateq}) is a matrix Riccatti equation and has general solution \cite{zwill}
\[
V=Q\left[V_{0}^{-1}-\int_{0}^{s}Q^{T}AQds\right]^{-1}Q^{T},
\]
where $V_0=V(0)$ and $Q$ is a $2\times2$ matrix satisfying the equation
\[
Q'=BQ, \qquad Q(0)=I,
\]
$I$ being the $2\times 2$ identity matrix. 

This has the solution
\[
Q=\begin{bmatrix} e^{\frac{\gamma}{2}s} & 0 \\ & \\ 0 & e^{\frac{\bar{\gamma}}{2}s} \end{bmatrix},
\qquad \mbox{and so} \qquad
Q^{T}AQ=\begin{bmatrix} \bar{\alpha} e^{\gamma s} & 0 \\ & \\ 0 & \alpha e^{\bar{\gamma}s} \end{bmatrix}.
\]
Thus
\[
\int_{0}^{s}Q^{T}AQds= \begin{bmatrix} \frac{\bar{\alpha}}{\gamma}( e^{\gamma s}-1) & 0 \\ & \\ 0 & \frac{\alpha}{\bar{\gamma}}(e^{\bar{\gamma}s}-1) \end{bmatrix}.
\]
If the initial value of $V^{-1}$ is 
\[
V_{0}^{-1}=\begin{bmatrix} v_{1} & v_{2} \\ -v_{2} & \bar{v}_{1} \end{bmatrix},
\]
for $v_{2}\in {\mathbb{R}}$ and $v_{1}\in {\mathbb{C}}$ then
 \[
C=V_{0}^{-1}-\int_{0}^{s}Q^{T}AQds=\begin{bmatrix}  v_{1}-\frac{\bar{\alpha}}{\gamma}( e^{\gamma s}-1) & v_{2}\\ & \\ -v_{2} & \bar{v}_{1}-\frac{\alpha}{\bar{\gamma}}(e^{\bar{\gamma}s}-1) \end{bmatrix},
\]
and the determinant of $C$ is
\[
\mbox{det}C=|v_{1}-\frac{\bar{\alpha}}{\gamma}( e^{\gamma s}-1)|^2+v_{2}^2.
\]
So, we obtain the inverse of $C$:
\[
\left[V_{0}^{-1}-\int_{0}^{s}Q^{T}AQds\right]^{-1}=\frac{1}{\mbox{det}C}\begin{bmatrix} \bar{v}_{1}-\frac{\alpha}{\bar{\gamma}}(e^{\bar{\gamma}s}-1)  & -v_{2}\\ & \\ v_{2} & v_{1}-\frac{\bar{\alpha}}{\gamma}( e^{\gamma s}-1)  \end{bmatrix}.
\]
Finally we find the solution
\[
V=QC^{-1}Q^{T}=\frac{1}{\mbox{det}C}\begin{bmatrix} e^{\gamma s}\left(\bar{v}_{1}-\frac{\alpha}{\bar{\gamma}}(e^{\bar{\gamma}s}-1)\right) & -v_{2}e^{\frac{\gamma+\bar{\gamma}}{2}s} \\ & \\ v_{2}e^{\frac{\gamma+\bar{\gamma}}{2}s} & e^{\bar{\gamma}s}\left(v_{1}-\frac{\bar{\alpha}}{\gamma}(e^{\gamma s}-1)\right)  \end{bmatrix} .
\] 
For $\beta=0$ then, the integral curves are
\[
t(s)=-\frac{v_{2}}{|v_{1}-\frac{\bar{\alpha}}{\gamma}( e^{\gamma s}-1)|^2+v_{2}^2}e^{\frac{\gamma+\bar{\gamma}}{2}s},
\]
\[
z(s)=\frac{\bar{v}_{1}-\frac{\alpha}{\bar{\gamma}}(e^{\bar{\gamma}s}-1)}{|v_{1}-\frac{\bar{\alpha}}{\gamma}( e^{\gamma s}-1)|^2+v_{2}^2}e^{\gamma s}.
\]

To solve the general case with $\beta\neq 0$, choose a complex number $\tau$ satisfying 
$\beta=-\gamma\tau-\bar{\alpha}\tau^2$. Then a shift $z\rightarrow z+\tau$ and completing the squares on the
$z$ term of right-hand side of the second equation of (\ref{e:lieH}) reduces the equations to the system
(\ref{e:mateq}). Thus the general solution turns out to be
\[
t(s)=-\frac{v_{2}}{|v_{1}-\frac{\bar{\alpha}}{\gamma_{1}}( e^{\gamma_{1}s}-1)|^2+v_{2}^2}e^{\frac{\gamma_{1}+\bar{\gamma}_{1}}{2}s},
\]
\[
z(s)=\frac{\bar{v}_{1}-\frac{\alpha}{\bar{\gamma}_{1}}(e^{\bar{\gamma}_{1}s}-1)}{|v_{1}-\frac{\bar{\alpha}}{\gamma_{1}}( e^{\gamma_{1}s}-1)|^2+v_{2}^2}e^{\gamma_{1}s}+\tau ,
\]
where $\gamma_{1}=\gamma+2\bar{\alpha}\tau$.

Setting
\[
v_{1}=\frac{\bar{z}_{0}-\bar{\tau}}{|z_{0}-\tau|^2+t_{0}^2}\, , \quad v_{2}=-\frac{t_{0}}{|z_{0}-\tau|^2+t_{0}^2},
\]
we obtain that $z(0)=z_{0}$ and $t(0)=t_{0}$ and the result follows.

\end{pf}

Since the above transformations (($t_0,z_0$),$s$)$\mapsto$($t(s),z(s)$) are isometries of ${\mathbb{H}}^3$,
 they map oriented geodesics to 
oriented geodesics. The following describes explicitly this as a map from ${\mathbb{L}}({\mathbb{H}}^3)$ to itself.

\begin{Prop}\label{p:isomgpLH}
The above action maps oriented geodesics in $(\xi,\eta)$ coordinates according to: for $\alpha$ and $\gamma$ not both zero:
\[
\xi\mapsto \xi e^{-\bar{\gamma}_{1}s}\left[ \frac{\alpha}{\bar{\gamma}_{1}}(e^{\bar{\gamma}_{1}s}-1)\left(\bar{\eta}+\frac{1}{\xi}-\bar{\tau}\right)-1\right]\left[ \frac{\alpha}{\bar{\gamma}_{1}}(e^{\bar{\gamma}_{1}s}-1)\left(\bar{\eta}-\frac{1}{\xi}-\bar{\tau}\right)-1\right],
\]
\[
\eta\mapsto \frac{\eta-\tau-\frac{\bar{\alpha}}{\gamma_{1}}(e^{\gamma_{1}s}-1)\left(\eta+\frac{1}{\bar{\xi}}-\tau\right)\left(\eta-\frac{1}{\bar{\xi}}-\tau\right)}{\left[\frac{\bar{\alpha}}{\gamma_{1}}(e^{\gamma_{1}s}-1)\left(\eta+\frac{1}{\bar{\xi}}-\tau\right)-1\right]\left[\frac{\bar{\alpha}}{\gamma_{1}}(e^{\gamma_{1}s}-1)\left(\eta-\frac{1}{\bar{\xi}}-\tau\right)-1\right]}e^{\gamma_{1}s}+\tau ,
\]
and for $\alpha=\gamma=0$
\[
\xi\mapsto \xi \qquad \qquad \eta\mapsto \eta+\beta s.
\]
\end{Prop}
\begin{pf}
It suffices to work on ${\mathbb{U}}\subset{\mathbb{L}}({\mathbb{H}}^3)$ and consider the action on the oriented geodesics
\[
z_{0}=\eta + \frac{\tanh r}{\bar{\xi}}, \quad t_{0}=\frac{1}{\sqrt{\xi \bar{\xi}}\cosh r}.
\]
We will find the oriented geodesic that is obtained by mapping this oriented geodesic by the 1-parameter group of 
actions in the last Proposition. 

The case $\alpha=\gamma=0$ follows trivially, so we omit the proof and 
consider the case where $\alpha$ and $\gamma$ are not both zero.

Denote 
\begin{align}
\Lambda&=|z_{0}-\tau|^2+t_{0}^2=\left(\eta-\tau +\frac{\tanh r}{\bar{\xi}}\right)\left(\bar{\eta}-\bar{\tau} +\frac{\tanh r}{\xi}\right),\nonumber \\
&=|\eta-\tau|^2+\frac{1}{|\xi|^2}+\left(\frac{\eta-\tau}{\xi}+\frac{\bar{\eta}-\bar{\tau}}{\bar{\xi}}\right)\tanh r.\nonumber
\end{align}
Then we obtain $v_{1}\bar{v}_{1}+v_{2}^2=\Lambda^{-1}$. Set $\lambda=\alpha\gamma_{1}^{-1}(e^{\gamma_{1}s}-1)$ and
\begin{align}
\left|v_{1}-\lambda\right|^2+v_{2}^2=&\frac{1-\lambda (z_{0}-\tau)-\bar{\lambda} (\bar{z}_{0}-\bar{\tau})+|\lambda|^2\Lambda}{\Lambda},\nonumber \\
=&\frac{1-\lambda(\eta-\tau)-\bar{\lambda}(\bar{\eta}-\bar{\tau})+|\lambda|^2\left(|\eta-\tau|^2+\frac{1}{|\xi|^2}\right)}
   {\Lambda},\nonumber\\
&\qquad\qquad+\frac{\left(\frac{|\lambda|^2(\bar{\eta}-\bar{\tau})-\lambda}{\bar{\xi}}+\frac{|\lambda|^2(\eta-\tau)-\bar{\lambda}}{\xi}\right)\tanh r}{\Lambda}.\nonumber
\end{align}
We now shift $\eta$ to $\eta+\tau$ - this will simplify the calculations and can be undone at the end by shifting
$\eta$ to $\eta-\tau$. Thus the above expression becomes
\[
|v_{1}-\lambda|^2+v_{2}^2=\frac{1-\lambda\eta-\bar{\lambda}\bar{\eta}+|\lambda|^2\left(|\eta|^2+\frac{1}{|\xi|^2}\right)
+\left(\frac{|\lambda|^2\bar{\eta}-\lambda}{\bar{\xi}}+\frac{|\lambda|^2\eta-\bar{\lambda}}{\xi}\right)\tanh r}{\Lambda}.
\]
Then,
\[
t=-\frac{v_{2}}{|v_{1}-\lambda|^2+v_{2}^2}e^{\frac{\gamma_{1}+\bar{\gamma}_{1}}{2}s},
\]
and therefore
\[
t=\frac{|\xi|^{-1}e^{\frac{\gamma_{1}+\bar{\gamma}_{1}}{2}s}}{\left[\left(1-\lambda\eta-\bar{\lambda}\bar{\eta}+|\lambda|^2|\eta|^2+\frac{|\lambda|^2}{|\xi|^2}\right)\cosh r+\left( \frac{|\lambda|^2\bar{\eta}-\lambda}{\bar{\xi}}+\frac{|\lambda|^2\eta-\bar{\lambda}}{\xi}  \right)\sinh r\right]}.
\]
Introducing values $A_{1},A_{2}$ by
\[
A_{1}=\left|\frac{\bar{\alpha}}{\gamma_1}(e^{\gamma_{1}s}-1)\left(\eta+\frac{1}{\bar{\xi}}\right)-1\right|, \quad A_{2}=\left|\frac{\bar{\alpha}}{\gamma_1}(e^{\gamma_{1}s}-1)\left(\eta-\frac{1}{\bar{\xi}}\right)-1\right|,
\]
then we may rewrite $t$ as follow
\[
t=\frac{1}{\sqrt{(\xi e^{-\bar{\gamma}_{1}s})(\bar{\xi}e^{-\gamma_{1}s})A_{1}^2 A_{2}^2}\cosh \left(r+\log\frac{A_{1}}{A_{2}}\right)}.
\]
Hence, we have the following map
\begin{align}
\xi&\mapsto \xi e^{-\bar{\gamma}_{1}s}\left[ \frac{\alpha}{\bar{\gamma}_1}(e^{\bar{\gamma}_{1}s}-1)\left(\bar{\eta}+\frac{1}{\xi}\right)-1\right]\left[ \frac{\alpha}{\bar{\gamma}_1}(e^{\bar{\gamma}_{1}s}-1)\left(\bar{\eta}-\frac{1}{\xi}\right)-1\right],\nonumber \\
r&\mapsto r+\log\left( \frac{A_{1}}{A_{2}}\right).\nonumber
\end{align}
It remains now to find where does the $\eta$ is mapped. To find this we have to find the $z$.
\begin{align}
z&=\frac{\bar{v}_{1}-\frac{\alpha}{\bar{\gamma}_1}(e^{\bar{\gamma}_{1}s}-1)}{|v_{1}-\frac{\bar{\alpha}}{\gamma_1}( e^{\gamma_{1}s}-1)|^2+v_{2}^2}e^{\gamma_{1}s},\nonumber \\
&=\frac{z_{0}-\frac{\alpha}{\bar{\gamma}_1}(e^{\bar{\gamma}_{1}s}-1)\Lambda}{1-\lambda\eta-\bar{\lambda}\bar{\eta}+|\lambda|^2 |\eta|^2+\frac{|\lambda|^2 }{|\xi|^2}+\left(\frac{|\lambda|^2\bar{\eta}-\lambda}{\bar{\xi}}+\frac{|\lambda|^2\eta-\bar{\lambda}}{\xi}\right)\tanh r}e^{\gamma_{1}s},\nonumber \\
&=\frac{\left(\eta+\frac{\tanh r}{\bar{\xi}}-\frac{\alpha}{\bar{\gamma}_1}(e^{\bar{\gamma}_{1}s}-1)\left[|\eta|^2+\frac{1}{|\xi|^2}+\left(\frac{\eta}{\xi}+\frac{\bar{\eta}}{\bar{\xi}}\right)\tanh r\right]\right)e^{\gamma_{1}s}\cosh r}{\left(1-\lambda\eta-\bar{\lambda}\bar{\eta}+|\lambda|^2 |\eta|^2+\frac{|\lambda|^2 }{|\xi|^2}\right)\cosh r+\left(\frac{|\lambda|^2\bar{\eta}-\lambda}{\bar{\xi}}+\frac{|\lambda|^2\eta-\bar{\lambda}}{\xi}\right)\sinh r} ,\nonumber \\
&=\frac{\eta-\frac{\bar{\alpha}}{\gamma_1}(e^{\gamma_{1}s}-1)\left(\eta+\frac{1}{\bar{\xi}}\right)\left(\eta-\frac{1}{\bar{\xi}}\right)}{\left[\frac{\bar{\alpha}}{\gamma_1}(e^{\gamma_{1}s}-1)\left(\eta+\frac{1}{\bar{\xi}}\right)-1\right]\left[\frac{\bar{\alpha}}{\gamma_1}(e^{\gamma_{1}s}-1)\left(\eta-\frac{1}{\bar{\xi}}\right)-1\right]}e^{\gamma_{1}s}\nonumber \\
&\qquad\qquad+\frac{\tanh \left(r+\log \frac{A_{1}}{A_{2}}\right)}{\bar{\xi} e^{-\gamma_{1}s}\left( \frac{\bar{\alpha}}{\gamma_1}(e^{\gamma_{1}s}-1)\left(\eta+\frac{1}{\bar{\xi}}\right)-1\right)\left( \frac{\bar{\alpha}}{\gamma_1}(e^{\gamma_{1}s}-1)\left(\eta-\frac{1}{\bar{\xi}}\right)-1\right)},\nonumber
\end{align}
Therefore we have the mapping 
\begin{align}
\eta&\longmapsto \frac{\eta-\frac{\bar{\alpha}}{\gamma_1}(e^{\gamma_{1}s}-1)\left(\eta+\frac{1}{\bar{\xi}}\right)\left(\eta-\frac{1}{\bar{\xi}}\right)}{\left[\frac{\bar{\alpha}}{\gamma_1}(e^{\gamma_{1}s}-1)\left(\eta+\frac{1}{\bar{\xi}}\right)-1\right]\left[\frac{\bar{\alpha}}{\gamma_1}(e^{\gamma_{1}s}-1)\left(\eta-\frac{1}{\bar{\xi}}\right)-1\right]}e^{\gamma_{1}s},\nonumber \\
\xi&\longmapsto \xi e^{-\bar{\gamma}_{1}s}\left( \frac{\alpha}{\bar{\gamma}_1}(e^{\bar{\gamma}_{1}s}-1)\left(\bar{\eta}+\frac{1}{\xi}\right)-1\right)\left( \frac{\alpha}{\bar{\gamma}_1}(e^{\bar{\gamma}_{1}s}-1)\left(\bar{\eta}-\frac{1}{\xi}\right)-1\right),\nonumber \\
r&\longmapsto r+\log\left( \frac{A_{1}}{A_{2}}\right).\nonumber
\end{align}
Finally we shift $\eta$ back to $\eta-\tau$ to yield the stated result:
\[
\eta\longmapsto \frac{\eta-\tau-\frac{\bar{\alpha}}{\gamma_1}(e^{\gamma_{1}s}-1)\left(\eta+\frac{1}{\bar{\xi}}-\tau\right)\left(\eta-\frac{1}{\bar{\xi}}-\tau\right)}{\left[\frac{\bar{\alpha}}{\gamma_1}(e^{\gamma_{1}s}-1)\left(\eta+\frac{1}{\bar{\xi}}-\tau\right)-1\right]\left[\frac{\bar{\alpha}}{\gamma_1}(e^{\gamma_{1}s}-1)\left(\eta-\frac{1}{\bar{\xi}}-\tau\right)-1\right]}e^{\gamma_{1}s}+\tau ,
\]
\[
\xi\longmapsto \xi e^{-\bar{\gamma}_{1}s}\left[ \frac{\alpha}{\bar{\gamma}_1}(e^{\bar{\gamma}_{1}s}-1)\left(\bar{\eta}+\frac{1}{\xi}-\bar{\tau}\right)-1\right]\left[ \frac{\alpha}{\bar{\gamma}_1}(e^{\bar{\gamma}_{1}s}-1)\left(\bar{\eta}-\frac{1}{\xi}-\bar{\tau}\right)-1\right].
\]
\end{pf}

\noindent{\bf Proof of Theorem 3}:
We have obtained a map 
\[
F\colon {\mathbb{L}}({\mathbb{H}}^3)\times {\mathbb{R}}\rightarrow {\mathbb{L}}({\mathbb{H}}^3)\colon ((\xi,\eta),s)\longmapsto (f_{1}(\xi,\eta,s),f_{2}(\xi,\eta,s)).
\]
Where, for $\alpha$ and $\gamma$ not both zero:
\begin{align}
f_{1}(\xi,\eta,s)&=\xi e^{-\bar{\gamma}_{1}s}\left[ \frac{\alpha}{\bar{\gamma}_1}(e^{\bar{\gamma}_{1}s}-1)\left(\bar{\eta}+\frac{1}{\xi}-\bar{\tau}\right)-1\right]\left[ \frac{\alpha}{\bar{\gamma}_1}(e^{\bar{\gamma}_{1}s}-1)\left(\bar{\eta}-\frac{1}{\xi}-\bar{\tau}\right)-1\right],\nonumber \\
f_{2}(\xi,\eta,s)&=\frac{\eta-\tau-\frac{\bar{\alpha}}{\gamma_1}(e^{\gamma_{1}s}-1)\left(\eta+\frac{1}{\bar{\xi}}-\tau\right)\left(\eta-\frac{1}{\bar{\xi}}-\tau\right)}{\left[\frac{\bar{\alpha}}{\gamma_1}(e^{\gamma_{1}s}-1)\left(\eta+\frac{1}{\bar{\xi}}-\tau\right)-1\right]\left[\frac{\bar{\alpha}}{\gamma_1}(e^{\gamma_{1}s}-1)\left(\eta-\frac{1}{\bar{\xi}}-\tau\right)-1\right]}e^{\gamma_{1}s}+\tau, \nonumber 
\end{align}
and for $\alpha=\gamma=0$
\[
f_1(\xi,\eta,s)=\xi\qquad\quad f_2(\xi,\eta,s)=\eta+\beta s.
\]

We will show that for any $s\in \mathbb{R}$ this map is an isometry of ${\mathbb{G}}$. In order to do so, 
we first find the derivative of $F$ at the point $s=0$
\[
F_{*}\left(\left.\frac{d}{ds}\right|_0\right)=\left.\frac{df_{1}}{ds}\right|_0\frac{\partial}{\partial \xi}+\left.\frac{d\bar{f}_{1}}{ds}\right|_0\frac{\partial}{\partial \bar{\xi}}+\left.\frac{df_{2}}{ds}\right|_0\frac{\partial}{\partial \eta}+\left.\frac{d\bar{f}_{2}}{ds}\right|_0\frac{\partial}{\partial \bar{\eta}}.
\]

Calculating the derivatives of $F$ in the case $\alpha=\gamma=0$ and converting to ($\mu_1,\mu_2$) coordinates we obtain 
the vector field
\[
{\mathbb{K}}=F_{*}\left(\left.\frac{d}{ds}\right|_0\right)={\mathbb{R}}\mbox{e}\left(-\beta\frac{\partial}{\partial \mu_{1}}-\bar{\beta}\mu_2^2\frac{\partial}{\partial \mu_{2}}\right),
\]
which are the Killing vectors of ${\mathbb{G}}$ found in Proposition \ref{p:lieLH} with $c_1=-\beta$, 
$c_2=0$ and $c_3=0$.

On the other hand, suppose that $\alpha$ and $\gamma$ are not both zero. Then
\[
\left.\frac{df_{1}}{ds}\right|_0=-\bar{\gamma}_{1}\xi-2\alpha\xi(\bar{\eta}-\bar{\bar{\tau}}),\qquad \left.\frac{df_{2}}{ds}\right|_0=\gamma_{1}(\eta-\tau)+\bar{\alpha}\left[(\eta-\tau)^2+\frac{1}{\bar{\xi}^2}\right].
\]
Recalling now the change of coordinates to $(\mu_{1},\mu_{2})$ in ${\mathbb{L}}({\mathbb{H}}^3)$ we obtain
\begin{align}
{\mathbb{K}}^{\mu_{1}}&=-\frac{1}{4}\left(\mu_{1}+\frac{1}{\bar{\mu}_{2}}\right)^2\left.\frac{d\bar{f}_{1}}{ds}\right|_0-\left.\frac{df_{2}}{ds}\right|_0=\gamma_{1}\tau-\bar{\alpha}\tau^2+(\gamma_{1}-2\bar{\alpha}\tau)\mu_{1}-\bar{\alpha}\mu_{1}^2,\nonumber \\
{\mathbb{K}}^{\mu_{2}}&=\frac{\mu_{2}^2}{4}\left(\bar{\mu}_{1}+\frac{1}{\mu_{2}}\right)^2\left.\frac{df_{1}}{ds}\right|_0-\mu_{2}^2\left.\frac{d\bar{f}_{2}}{ds}\right|_0=-\alpha-(\bar{\gamma}_{1}-2\alpha\bar{\tau})\mu_{2}+\left(\bar{\gamma}_{1}\bar{\tau}-\alpha\bar{\tau}^2\right)\mu_{2}^2.\nonumber
\end{align}

Now substituting $\gamma_1=\gamma+2\bar{\alpha}\tau$ and $\gamma\tau+\bar{\alpha}\tau^2=-\beta$ back in, we find that
\[
{\mathbb{K}}^{\mu_{1}}=-\beta+\gamma\mu_{1}-\bar{\alpha}\mu_{1}^2,\qquad\qquad
{\mathbb{K}}^{\mu_{2}}=-\alpha-\bar{\gamma}\mu_{2}-\bar{\beta}\mu_{2}^2.
\]
These are precisely the Killing vectors of ${\mathbb{G}}$ found in Proposition \ref{p:lieLH} with $c_1=-\beta$, 
$c_2=\gamma$ and $c_3=-\bar{\alpha}$, and hence
\[
\mbox{Iso}_0({\mathbb{L}}({\mathbb{H}}^3),{\mathbb{G}})=\mbox{Iso}_0({\mathbb{H}}^3,g),
\]
as claimed.

\hspace{4.0in}$\Box$

\section{The Geodesics of ${\mathbb{G}}$}

A curve in ${\mathbb{L}}({\mathbb{H}}^3)$ is a 1-parameter family of oriented geodesics in ${\mathbb{H}}^3$ - which we 
refer to as a {\it ruled surface} in ${\mathbb{H}}^3$. The ruled surfaces that come from geodesics of the neutral
K\"ahler metric in ${\mathbb{L}}({\mathbb{H}}^3)$ have a particularly elegant characterisation:

\begin{Thm}
The geodesics of the K\"ahler metric ${\mathbb{G}}$ are generated by the 1-parameter subgroups of the isometry group
of  ${\mathbb{G}}$. 

A ruled surface generated by a geodesic of ${\mathbb{G}}$ is a minimal surface in ${\mathbb{H}}^3$, and the
geodesic is null iff the ruled surface is totally geodesic.
\end{Thm}
\begin{pf}
For $I\subset{\mathbb{R}}$, let $c: I\rightarrow {\mathbb{L}}({\Bbb{H}}^3)$ be a geodesic in ${\mathbb{L}}({\Bbb{H}}^3)$ 
with affine parameter $t$. By an isometry we can move the geodesic to lie in ${\mathbb{U}}\subset{\mathbb{L}}({\Bbb{H}}^3)$
and use coordinates ($\xi,\eta$) as earlier. Thus the geodesic is given by $c(t)=(\xi(t),\eta(t))$ satisfying the geodesic 
equations
\[
\nabla_{\dot{c}}\;\dot{c}=0,
\]
where $\dot{c}$ is the tangent vector to $c$:
\[
\dot{c}=\dot{\xi}(t)\frac{\partial}{\partial \xi}+\dot{\eta}(t)\frac{\partial}{\partial \eta}+\dot{\overline{\xi}}(t)\frac{\partial}{\partial \overline{\xi}}+\dot{\overline{\eta}}(t)\frac{\partial}{\partial \overline{\eta}},
\]
the dot denoting differentiation with respect of $t$.

These equations, using the metric expression in Proposition \ref{p:kaehlermetric} turn out to be
\begin{align}\label{e:geou}
\xi\ddot{\xi}-\dot{\xi}^2+(\dot{\overline{\eta}})^2 \xi^4=0, \qquad \overline{\xi}\ddot{\eta}+2\dot{\overline{\xi}}\dot{\eta}=0,
\end{align}
which we solve as follows. From the second equation we obtain
\[
\frac{d}{dt}\left(\overline{\xi}^2 \frac{d\eta}{dt}\right)=0,
\]
from which we get that
\[
\dot{\eta}=\frac{b_{3}}{\overline{\xi}^2},
\]
for some complex constant $b_3$.

Substituting this into the first equation of (\ref{e:geou}) we obtain 
\[
\xi\ddot{\xi}-\dot{\xi}^2+\overline{b}_{3}^2=0.
\]
This has solution
\[
\xi(t)=\frac{\overline{b}_{3}\sinh(b_{2}t+b_{1})}{b_{2}}.
\]
We now get that
\[
\dot{\eta}=\frac{\overline{b}_{2}^2}{b_{3}\sinh^2(\overline{b}_{2}t+\overline{b}_{1})},
\]
which we integrate to find the geodesics:

\[
\xi(t)=\frac{\overline{b}_{3}\sinh(b_{2}t+b_{1})}{b_{2}},\qquad \eta(t)=b_{4}-\frac{\overline{b}_{2}\cosh(\overline{b}_{2}t+\overline{b}_{1})}{b_{3}\sinh(\overline{b}_{2}t+\overline{b}_{1})},
\]
for $b_{1},b_{2},c_{3},c_{4}\in{\mathbb{C}}$. We note that the length of the tangent vector to the geodesic is
the constant $i(\bar{b}^2_2-b_2^2)/4$.

The geodesic in terms of ($\mu_{1},\mu_{2}$) coordinates is given by
\[
\mu_1=-b_4+\frac{\overline{b}_{2}[1+\cosh(\overline{b}_{2}t+\overline{b}_{1})]}{b_3 \sinh (\overline{b}_{2}t+\overline{b}_{1})},\quad \mu_{2}=\frac{\overline{b}_{3}\sinh(b_{2}t+b_{1})}{\overline{b}_{3}\overline{b}_{4}\sinh(b_{2}t+b_{1})+b_{2}[1-\cosh(b_{2}t+b_{1})]}.
\]
The tangent vector of the geodesic is
\[
\dot{c}(t)={\mathbb{R}}\mbox{e}\left(\dot{\mu}_{1}\frac{\partial}{\partial \mu_{1}}
    +\dot{\mu}_{2}\frac{\partial}{\partial \mu_{2}}\right),
\]
where,
\[
\dot{\mu}_{1}=-\frac{\overline{b}_{2}^2 [1+\cosh(\overline{b}_{2}t+\overline{b}_{1})]}{b_{3}\sinh^2(\overline{b}_{2}t+\overline{b}_{1})},\quad  \dot{\mu}_{2}=\frac{b_{2}^2 \overline{b}_{3}[\cosh(b_{2}t+b_{1})-1]}{\left(\overline{b}_{3}\overline{b}_{4}\sinh(b_{2}t+b_{1})+b_{2}[1-\cosh(b_{2}t+b_{1})]\right)^2}.
\]
Thus for all $t$
\[
\dot{\mu}_{1}(t)=c_{1}+c_{2}\mu_{1}(t)+c_{3}\mu_{1}^2(t),\qquad \dot{\mu}_{2}(t)=\overline{c}_{3}-\overline{c}_{2}\mu_{2}(t)+\overline{c}_{1}\mu_{2}^2(t),
\]
where,
\[
c_{1}=-\frac{b_{3}^2 b_{4}^2-\overline{b}_{2}^2}{2b_{3}},\qquad c_{2}=-b_{3}b_{4},\qquad c_{3}=-\frac{b_{3}}{2}.
\]
We conclude from Proposition \ref{p:lieLH}, that the tangent vector of the geodesic is the restricion of a killing vector 
of ${\mathbb{G}}$. 

We now compute the second fundamental form of the ruled surface in ${\mathbb{H}}^3$ generated by the geodesics
of ${\mathbb{G}}$. In order to simplify the calculations we first utilise an isometry of ${\mathbb{U}}$ ({\it cf}. 
Proposition \ref{p:isomgpLH}) with
\[
\alpha=0,\qquad s=1, \qquad e^{\gamma_{1}}=b_{3},\qquad {\mbox{and}} \qquad\tau=\frac{b_{3}b_{4}}{b_{3}-1}.
\]
This simplifies the geodesic to
\[
\xi(t)=\frac{\sinh(b_{2}t+b_{1})}{b_{2}},\qquad \eta(t)=-\frac{\overline{b}_{2}\cosh(\overline{b}_{2}t+\overline{b}_{1})}
  {\sinh(\overline{b}_{2}t+\overline{b}_{1})}.
\]
The map $\Phi\colon {\mathbb{U}}\times {\mathbb{R}}\rightarrow {\mathbb{H}}^3\colon ((\xi,\eta),r)\mapsto (x^0,x^1,x^2)$, 
where 
\[
x^0=\frac{1}{|\xi|\cosh r},\qquad x^1=\frac{\eta+\overline{\eta}}{2}+\frac{\tanh r}{2}\left(\frac{1}{\xi}+\frac{1}{\overline{\xi}}\right),\]
\[
x^2=\frac{i(\overline{\eta}-\eta)}{2}+\frac{i\tanh r}{2}\left(\frac{1}{\xi}-\frac{1}{\overline{\xi}}\right),
\] 
now yields the parameterisation of the ruled surface S in ${\mathbb{H}}^3$. In particular, we have
a surface given by $(r,t)\mapsto (x^{0}(r,t),x^{1}(r,t),x^{2}(r,t))$ which has induced metric
\[
\overline{g}_{lk}=g_{ij}\frac{\partial x^{i}}{\partial y^{l}}\frac{\partial x^{j}}{\partial y^{k}},
\] 
and normal vector $\overline{N}$
\begin{align}
\overline{N}=&\left(\frac{\partial x^1}{\partial y^1}\frac{\partial x^2}{\partial y^2}-\frac{\partial x^1}{\partial y^2}\frac{\partial x^2}{\partial y^1}\right)\frac{\partial}{\partial x^0}-\left(\frac{\partial x^0}{\partial y^1}\frac{\partial x^2}{\partial y^2}-\frac{\partial x^0}{\partial y^2}\frac{\partial x^2}{\partial y^1}\right)\frac{\partial}{\partial x^1} \nonumber \\
&\qquad\qquad+\left(\frac{\partial x^0}{\partial y^1}\frac{\partial x^1}{\partial y^2}-\frac{\partial x^0}{\partial y^2}\frac{\partial x^1}{\partial y^1}\right)\frac{\partial}{\partial x^2},\nonumber
\end{align}
where $y^1=r$ and $y^2=t$.

The second fundamental form $K_{ab}$ of S is given by
\[
K_{ab}=\frac{\partial x^{k}}{\partial y^{a}}\left(\frac{\partial N_{k}}{\partial y^{b}}-\tilde{\Gamma}^{i}_{lk}\frac{\partial x^{l}}{\partial y^{b}}N_{i}\right),\qquad a,b=1,2,
\]
where $N$ is the unit normal of S and $\tilde{\Gamma}^{i}_{lk}$ are the christoffel symbols of the
hyperbolic metric $g$. 

After some lengthy calculations, the components of the second fundamental form of the ruled surface S in ${\mathbb{H}}^3$ 
are found to be:
\[
K_{r r}=0, \qquad K_{r t}= \frac{|\sinh(b_{2}t+b_{1})|}{2M}(b^{2}_{2}-\overline{b}_{2}^{2}),
\]
\[
K_{t t}=\frac{{\mathbb{R}}\mbox{e}[\overline{b}_{2}\sinh(b_{2}t+b_{1})]}{|\sinh(b_{2}t+b_{1})|M}(b^{2}_{2}-\overline{b}_{2}^{2}),
\]
where
\begin{align}
M^2=&2|b_{2}|^2 \cosh r \sinh r \; {\mathbb{R}}\mbox{e}[\cosh (b_{2}t+b_{1})]-4|b_{2}|^2\sinh^2 r\; 
     (|\cosh (b_{2}t+b_{1})|^2+1)   \nonumber\\
&\qquad\qquad-2|b_{2}|^2|\cosh (b_{2}t+b_{1})|^2-|\sinh(b_{2}t+b_{1})|^2\; (b_{2}^2+\overline{b}_{2}^2)-2|b_{2}|^2,\nonumber
\end{align}
Computing the mean curvature of S we get
\[
H=\overline{g}^{ab}K_{ab}=0
\]
and hence every geodesic in ${\mathbb{L}}({\mathbb{H}}^3)$ is a ruled minimal surface in  ${\mathbb{H}}^3$.

In addition, we see that the second fundamental form vanishes when $b_{2}\in {\mathbb{R}}$ or $b_{2}\in i{\mathbb{R}}$.
As noted earlier, the length of the tangent vector to the geodesic is $i(\bar{b}^2_2-b_2^2)/4$, and so the geodesic
in ${\mathbb{L}}({\mathbb{H}}^3)$ is null iff the ruled surface S in ${\mathbb{H}}^3$ is totally geodesic. 
\end{pf}

\vspace{0.1in}
\setlength{\epsfxsize}{4.5in}
\begin{center}
   \mbox{\epsfbox{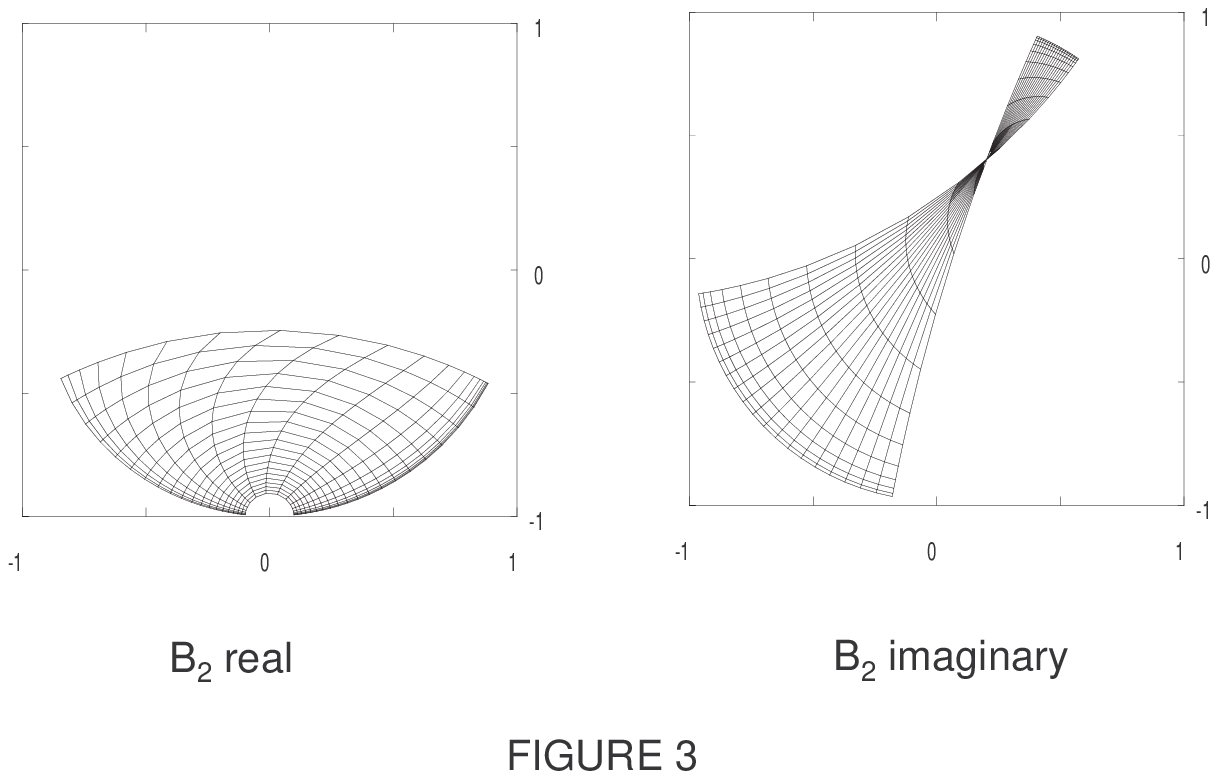}}
\end{center}
\vspace{0.1in}

In Figure 3 we show the ruled surfaces generated by two null geodesics in the ball model. These geodesics lie in the 
equatorial plane and 
correspond to a hyperbolic translation and a hyperbolic rotation. In Figure 4 we show a non-null geodesic, which
generates a hyperbolic helicoid.

\vspace{0.1in}
\setlength{\epsfxsize}{3.5in}
\begin{center}
   \mbox{\epsfbox{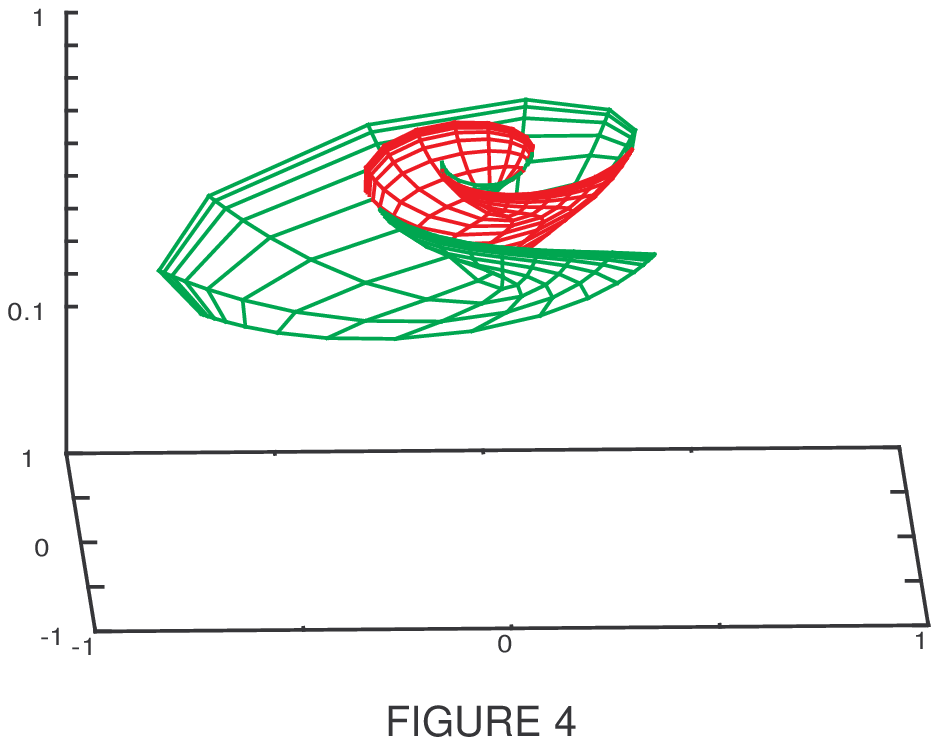}}
\end{center}


\vspace{0.1in}



\begin{thebibliography}{10}
\bibitem{gak4}
B. Guilfoyle and W. Klingenberg, {\it An indefinite K\"ahler metric on the space of oriented lines}, J. London Math.
Soc. {\bf 72} (2005), 497--509.
\bibitem{gak5}
B. Guilfoyle and W. Klingenberg, {\it A neutral K\"ahler metric on the space of time-like lines in Lorentzian 3-spaces} 
    (2006) math.DG/0608782.
\bibitem{gak6}
B. Guilfoyle and W. Klingenberg, {\it A neutral K\"ahler surface with applications in geometric optics}, in New Results
in Lorentz Geometry, World Scientific, Vienna (to appear).
\bibitem{hitch}
N.J. Hitchin, {\it Monopoles and geodesics}, Comm. Math. Phys. {\bf 83} (1982), 579--602.
\bibitem{kan}
S. Kobayashi and K. Nomizu, Foundations of Differential Geometry, Volume II, Wiley and Sons, New York, 1996.
\bibitem{kok}
M. Kokubu {\it et al}, {\it Singularities of flat fronts in hyperbolic space},  Pacific J. Math. {\bf 221} (2005), 303--352.\bibitem{Rat}
J. G. Ratcliffe, Foundations of hyperbolic manifolds, Springer-Verlag, New York, 1994.
\bibitem{roit}
P. Roitman, {\it Flat surfaces in hyperbolic space as normal surfaces to a congruence of geodesics}, Tohoku Math. J. 
(to appear).
\bibitem{salvai1}
M. Salvai, {\it On the geometry of the space of oriented lines in Euclidean space}, Manuscripta Math. {\bf 118} (2005),
181--189.
\bibitem{salvai2}
M. Salvai, {\it On the geometry of the space of oriented lines of hyperbolic space}, (2007)  math.DG/0702365.
\bibitem{small}
A. Small, {\it Surfaces of constant mean curvature 1 in ${\mathbb{H}}^3$ and algebraic curves on a quadric}, Proc. of the 
A.M.S. {\bf 122} (1994), 1211--1220.
\bibitem{weier}
K. Weierstrass, {\it Untersuchungen \"{u}ber die Fl\"{a}chen, deren mittlere Kr\"{u}mmung \"{u}berall gleich Null ist}, Monatsber. Akad. Wiss. Berlin (1866), 612-625. 
\bibitem{whitt}
E. T. Whittaker, {\it On the partial differential equations of mathematical physics }, Math. Ann. {\bf 57} (1903), 333-355.
\bibitem{zwill}
D. Zwillinger, Handbook of Differential Equations, Academic Press, New York, 1989.






\end{thebibliography}
\end{document}